\documentclass{article}
\usepackage{graphicx, dsfont, amssymb, natbib, color} 

\usepackage{booktabs}
\usepackage{array}
\usepackage{multirow}
\usepackage{xcolor}
\usepackage{algorithm}
\usepackage{algorithmic}
\usepackage{url}
\usepackage{subcaption}
\usepackage[margin=1in]{geometry}
\usepackage{pdfpages}
\usepackage{amsbsy}
\usepackage{amssymb}
\usepackage{graphicx}
\usepackage{float}
\usepackage{graphicx} 
\usepackage{amsmath}
\usepackage{amsfonts}
\usepackage{hyperref}
\usepackage{authblk}
\usepackage{natbib}
\usepackage{texab}
\usepackage[utf8]{inputenc}
\usepackage{lmodern}
\usepackage{enumerate}
\DeclareUnicodeCharacter{2212}{-}
\usepackage{dsfont}
\usepackage{amsfonts}
\usepackage{amssymb}
\usepackage{mathtools}
\usepackage{relsize}
\usepackage{tikz}
\usetikzlibrary{arrows.meta,positioning,calc}
\usepackage{color} 
\usepackage{listings}
\definecolor{Mygrey}{gray}{0.75}
\definecolor{Cgrey}{gray}{0.4}
\newtheorem{theorem}{Theorem}
\newtheorem{lemma}{Lemma}
\newtheorem{proposition}{Proposition}

\newtheorem{definition}{Definition}

\newtheorem{remark}{Remark}

\title{Linear regression with known noise distribution up to a scale: The reward of not using the OLSE}
\author[1]{Fadoua Balabdaoui}
\author[2]{Justine Leclerc}

\affil[1]{Department of Mathematics, ETH Zürich, Zürich, Switzerland}
\begin{document}
\affil[2]{Center of Experimental Rheumatology, Department of Rheumatology, University Hospital Zurich, University of Zurich, Schlieren, Switzerland}

\maketitle

\begin{abstract}
While the ordinary least squares estimator (OLSE) is still the most widely used estimator in linear regression models, other estimators can be more efficient when the error distribution is not Gaussian. In this paper, our goal is to evaluate this efficiency in the case of the Maximum Likelihood estimator (MLE) when the noise distribution belongs to a scale family. Under some regularity conditions, we show that $(\widehat \beta_n, \widehat s_n)$, the MLE of the unknown regression vector $\beta_0$ and the scale parameter $s_0$, exists and derive the expression of the asymptotic efficiency of $\widehat \beta_n$ over the OLSE. For three given scale families of densities, we quantify the true statistical gain of the MLE as a function of their deviation from the Gaussian family. To illustrate the theory, we present simulation results for different settings and also compare the MLE to the OLSE for the real market fish dataset.

Keywords: linear regression, maximum likelihood estimation, non-Gaussian noise
\end{abstract}

\section{Introduction}
\subsection{Motivation and background} 
This article studies estimation methods for linear regression models when the noise distribution belongs to a scale family with unknown scale parameter. Our goal is to use this information to improve estimation performance. In particular, we seek estimators that are at least as efficient as the Ordinary Least Squares Estimator (OLSE), which does not use information about the noise distribution.
The present work is mainly motivated by the study of \citet{balabdaoui}. In that paper, the authors considered a linear regression setting where the response variables and predictor vectors are unmatched and may come from different sources. More precisely, they considered two independent datasets, $\{Y_i\}_{i=1}^{n_Y}$ and $\{X_i\}_{i=1}^{n_X}$, with possibly $n_Y \ne n_X$. As in \cite{balabdaoui}, we assume in the sequel that $n_X = n_Y$. In a linear framework, the two datasets are related through the equality in distribution
\begin{eqnarray}\label{unmatched}
Y \overset{d}{=}\beta_0^TX + \epsilon.
\end{eqnarray}
The authors showed that it is possible to estimate $\beta_0$ in the unmatched regression model (\ref{unmatched}) by solving a deconvolution problem. The resulting estimator is called the Deconvolution Least Squares Estimator (DLSE). However, consistency may fail because the model is not always identifiable. This happens when several values of $\beta_0$ lead to the same equality in distribution in (\ref{unmatched}).
For example, assume that $X \sim \mathcal N(0, \Sigma)$, where $\Sigma \in \mathbb R^{d \times d}$ is positive definite. Then the set of all $\beta_0$ satisfying (\ref{unmatched}) is the ellipsoid $\{ \beta \in \mathbb R^d:  \beta^T \Sigma \beta =  c \} $ for some constant $c > 0$. Hence, the true regression vector cannot always be uniquely identified.

To overcome this issue, \cite{balabdaoui} proposed a semi-supervised approach. The idea is to supplement the unmatched data with a smaller matched dataset $\displaystyle \{(Y_j,X_j)\}_{j=1}^{m}$. The term “matched” means that $X_j$ and $Y_j$ are observed for the same individual or item, as in the classical regression setting. The key assumption is that $\mathbb{E}[Y_j | X_j] =  \beta^T_0 X_j $ and that the error term $Y_j - \beta_0 X_j $ has the same distribution as $\epsilon$ in the unmatched model (\ref{unmatched}). 

\begin{figure}[h!]
\centering
\small
\begin{tikzpicture}[
    >=Stealth,
    box/.style={draw, rounded corners, text width=3.5cm, minimum height=1cm, align=center},
    smallbox/.style={draw, rounded corners, text width=3.5cm, minimum height=1cm, align=center},
    finalbox/.style={draw, rounded corners, text width=4.5cm, minimum height=1cm, align=center}
]

\node[box] (unmatched) {
\textbf{Unmatched dataset}\\
$(X_i)_{i=1}^n,\ (Y_i)_{i=1}^n$\\
};

\node[box, below=0cm of unmatched] (unmatched2) {
$Y_i \overset{d}{=} \beta_0^T X_i + \epsilon_i$\\
$F^\epsilon$ known\\
$\beta_0 \in \mathcal{B}_0$,\\
possibly $|\mathcal{B}_0|\geq 1$
};

\node[box, right=1.5cm of unmatched] (matched) {
\textbf{Matched dataset}\\
\smallskip
$(X_i,Y_i)_{i=1}^m$\\
};

\node[box, below=0cm of matched] (matched2) {
$Y_i = \beta_0^T X_i + \epsilon_i$\\
$F^\epsilon$ known
};

\node[box, left=1.5cm of unmatched] (semi) {
\textbf{Semi-supervised learning}\\
\smallskip
$n \gg m$
};

\node[smallbox, below=1cm of unmatched2] (dlse) {
\textbf{DLSE } $\widehat{\beta}_n$\\
\smallskip
(Azadkia and Balabdaoui)
};

\node[smallbox, below=1.5cm of matched2] (estimator) {
\textbf{Estimator } $\widehat{\beta}_m$ of $\beta_0$,\\
\smallskip
(MLE, alternative approach leveraging knowledge of $F^\epsilon$)
};

\node[finalbox, below=2.65cm of semi] (final) {
\textbf{New estimator } $\widehat{\beta}_{n,m}$ of $\beta_0$,\\
\smallskip
with good asymptotic properties
};


\draw[->, thick] (matched.west) to[bend right=25] (unmatched.east);

\draw[->, thick] (unmatched.west) -- (semi.east);

\draw[->, thick] (unmatched2) -- (dlse);
\draw[->, thick] (matched2) -- (estimator);
\draw[->, thick] (semi) -- (final);

\end{tikzpicture}
\caption{Semi-supervised learning design combining matched and unmatched datasets.}
\end{figure}

Even if the matched dataset is small, and the resulting estimator may therefore be less stable or less precise, it still provides useful information about the true parameter $\beta_0$. In particular, it helps identify the correct solution among the candidates obtained from the deconvolution step. This approach combines the strengths of both datasets: the matched data ensures identifiability, while the unmatched data can increase statistical power because of its potentially large sample size and the available knowledge about the noise distribution.

To derive the asymptotic properties of the DLSE, the authors of \citet{balabdaoui} assume that the noise distribution is known. However, they also mention that the method could be extended to cases where the noise distribution is either estimated from a sample of errors or assumed to belong to a scale family. In the latter case, the noise distribution is known up to an unknown scale parameter. As a result, the estimation procedure must recover both the regression coefficients and the scale parameter.

This naturally raises the question of whether one can construct an estimator that performs better than the OLSE in terms of asymptotic variance. Motivated by the large literature on asymptotic efficiency, a natural alternative to the OLSE is the Maximum Likelihood Estimator (MLE). Surprisingly, to the best of our knowledge, there are no results comparing the statistical efficiency of the MLE and the OLSE in this regression setting. In the next section, we review the literature on estimation in linear regression models with known noise distribution.

\subsection{Linear model with known error distribution: literature review } \label{littrev}

The least squares method only requires the existence of a second moment for the noise distribution. Under this assumption, the OLSE is the best linear unbiased estimator. When the noise is Gaussian, the OLSE is even optimal among all unbiased estimators of the regression vector. However, many studies have considered linear models where the noise distribution is known but not Gaussian. This naturally leads to the field of robust estimation, where Peter J. Huber (hereafter PJH) made major contributions.
In his 1964 paper \citep{huber1964}, PJH studied how sensitive the OLSE is to deviations from Gaussian assumptions. He was motivated by the poor performance of the OLSE under heavy-tailed or contaminated error distributions. To address this issue, he investigated alternative estimators that remain reliable under model misspecification.

As a first step, PJH focused on estimating a location parameter. Let $X_1, \ldots, X_n$ be i.i.d. observations from a distribution $F$. He introduced the notion of an \emph{M-estimator}, defined as,
\[
T_n = \arg\min_{t \in \mathbb{R}} \sum_{i=1}^n \rho(X_i - t),
\]
where $\rho: \mathbb{R} \to \mathbb{R}$ is a chosen loss function. This framework includes the sample mean when $\rho(u)=u^2$ and the sample median when $\rho(u)=|u|$. It also laid the foundation for a broad class of robust estimators.

In his 1972 monograph \citep{huber1972}, PJH extended the M-estimation framework from location models to linear regression. He studied the estimation of regression coefficients under model misspecification in the fixed design setting. Under suitable regularity conditions, he established consistency and asymptotic normality of the estimators.
This work was continued in 1973 \citep{huber1973}, where PJH derived explicit expressions for the asymptotic bias and variance of regression M-estimators. It is interesting to note that the variance estimator proposed there is similar to the one derived later in this article. Assuming that $\rho$ is twice continuously differentiable, the variance can be approximated by
\[
\frac{\frac{1}{n-d}\sum_{i=1}^n\left[\rho'\left(Y_i - \widehat{\beta}^\top X_i \right)\right]^2}
{\left(\frac{1}{n}\sum_{i=1}^n \rho''\left(Y_i - \widehat{\beta}^\top X_i \right)\right)^2},
\]
This estimator can be viewed as an early version of modern heteroskedasticity-consistent covariance estimators.
Finally, PJH’s 1981 book \citep{huber1981} summarizes many of his contributions to robust statistics, M-estimation, and maximum likelihood theory. Chapter 7 focuses on linear regression and studies the properties of the OLSE, including consistency and sensitivity to outliers. However, the framework considered there assumes fixed design matrices, unlike the random design setting studied in the present work.
Several other works have also studied inference in linear models with non-Gaussian errors.
For example, \citet{zeck} considered noise distributions from a power family parameterized by $(\theta,\mu,\sigma)\in(0,+\infty)\times\mathbb{R}\times(0,+\infty)$, with density
\begin{eqnarray*}
f^\epsilon(t)= k(\sigma, \theta)\exp\left(-\left|\frac{t-\mu}{\sigma}\right|^\theta\right), \qquad t\in\mathbb{R},
\end{eqnarray*}
where $k(\sigma, \theta)= \left(2\sigma \Gamma\left(1+\frac{1}{\theta}\right)\right)^{-1}$. The case $\theta=2$ corresponds to Gaussian noise, while smaller values of $\theta$ produce heavier tails. The authors mainly focused on methods that use the known form of the error distribution. They derived MLEs for both the regression coefficients and the parameter $\theta$, and compared them with the OLSE. Their results showed that the OLSE becomes less efficient and more biased as $\theta$ moves away from 2. They therefore argued that the MLE is preferable when the noise follows a power-family distribution.

In \citet{jlibene}, the authors studied a simple linear regression model with intercept, where the noise follows a uniform distribution $\mathcal{U}([-\theta,\theta])$ for some $\theta>0$. In the fixed design setting, the OLSE can be written as a linear combination of i.i.d. uniform random variables. The estimator therefore follows a generalized Irwin--Hall distribution. Assuming that $\theta$ is known, the authors derived the exact distribution of the OLSE and studied its asymptotic properties. They also proposed an unbiased estimator of the variance parameter. In particular, they showed that using the uniform noise assumption leads to more accurate confidence intervals than relying only on asymptotic normality.

Similarly, \citet{nolan} observed that the performance of the OLSE deteriorates under stable heavy-tailed error distributions. Through simulations, the authors showed that the OLSE is strongly affected by outliers, especially when the distribution differs substantially from the Gaussian case. By comparing root mean square deviations (RMSE), they found that stable MLE methods adapt much better to such settings. This makes the MLE a useful alternative when the noise distribution is highly dispersed.

Finally, in \citet{clancybecker}, the authors proposed an approximate MLE framework that accounts for uncertainty in both the design matrix and the error terms. In their setting, only the response variables are observed, while the predictor variables and the noise are known only through their distributions. Using moment generating functions and saddlepoint approximations, they approximated the likelihood function and recovered the regression vector. This framework is appealing because it goes beyond the classical fixed design assumption and allows for a more realistic regression model.

\subsection{Contributions and organization of the paper}

In the present article we address scenarios in linear regression where the design matrix is allowed to be random and the noise distribution is known up to a scale parameter. We consider the joint estimation of the unknown regression vector and the scale using the maximum likelihood approach.  Under suitable smoothness assumptions, we provide the following theoretical contributions:
\begin{enumerate}

\item We prove existence of the MLE of the regression vector and the scale parameter. While this step is most often overlooked or taken for granted, we give here a detailed argument that warrants that the likelihood admits at least a feasible maximizer.

\item We show consistency of any MLE (note that uniqueness does not necessarily hold). In the proof, we employ  techniques from empirical process theory. 

\item We prove the joint asymptotic normality of any MLE under some specified regularity assumptions on the density of the noise. Such assumptions allow us to exhibit the usual connection to the inverse of the Fisher Information matrix and also find that the estimators of the regression vector and the scale are asymptotically independent.

\item We quantify the asymptotic efficiency of the MLE of the regression vector when compared to the OLSE for three different parametric families. We illustrate that the statistical gain for using the MLE can be quite significant for a certain range of parameters.

\end{enumerate}

The manuscript is organized as follows. In Section \ref{sec: MLE}, we present the ML approach, and prove existence, consistency and asymptotic normality of any maximizer of the log-likelihood under specified regularity assumptions. 
In Section \ref{sec: AsympEffic} we show that the MLE of the regression vector is at least as efficient as the OLSE. Furthermore, we quantify this efficiency in case the noise distribution belongs to one of three scale families including the families of symmetric Gamma and Weibull distributions. To this aim, we exhibit numerically the parameter at which the MLE's efficiency is maximal. In Section \ref{simu}, we present simulations results in various settings as well as an application to the real fish market dataset to illustrate the theory. For a more fluid reading, all proofs are deferred to Appendix B.

\section{ML estimation with a scale-distribution for the noise} \label{sec: MLE}

\subsection{Context and assumptions} \label{context}
In this subsection, we outline the initial assumptions that will be maintained throughout the paper, except where explicitly stated otherwise.

Let $d \ge 1$ be an integer. We consider a random design where we observe $n$ independent pairs $(X_i, Y_i) \in \mathbb R^d \times \mathbb R$ for $i =1, \ldots, n$ the  (matched) linear regression model 
\begin{eqnarray}\label{LM}
Y = \beta^T_0  X + \epsilon
\end{eqnarray}
for some unknown regression vector $\beta_0 \in \mathbb R^d$. We will assume in the sequel that 
\begin{description}
        \item \textbf{(A1)} $X$ and $\epsilon$ are independent,
        \item \textbf{(A2)} the distribution of $X$ admits a density $f^X$ with respect to Lebesgue measure,   
        \item \textbf{(A3)}  the distribution of $\epsilon$ admits a density  $f^\epsilon$ with respect to Lebesgue measure which belongs to the scale family
        \begin{eqnarray}\label{scalefam}
            \left \{ \frac{1}{s} f\left(\frac{\cdot}{s}\right), s \in (0, \infty)    \right \}
        \end{eqnarray}
        where $f$ is some known even,  bounded and  continuous density such that $\int_{\mathbb R} t^2 f(t) dt < \infty$. 
        
\end{description}        

We will denote by $s_0$ the true but unknown scale parameter. In other words,  $f^\epsilon(t) = s^{-1}_0  f(t/s_0), t \in \mathbb R$.
In the context of \lq\lq classical\rq\rq \ linear regression, the noise distribution is typically assumed to be symmetric around 0, and to admit a finite variance.  In this regard, Assumption (A3) is not restrictive.  In the sequel, we will sometimes use the more convenient notation 
$f_s$ to mean $1/s f(\cdot / s)$ for $s > 0$ and $f_0$  to mean $f_{s_0}$, the true density of $\epsilon$.

\subsection{Existence}
 Using Assumptions (A1) and (A2), the joint density of $(X, Y)$ satisfying the linear model in (\ref{LM}) can be written as
$$
f^{(X, Y)}(x, y) =  f^X(x)  f^\epsilon(y - \beta^T_0  x)  = f^X(x) \frac{1}{s_0} f\left(\frac{y-\beta^T_0 x}{s_0}\right)
$$
for $(x, y)  \in \mathbb R^d \times \mathbb R$. Thus, based on the i.i.d. sample $(X_i, Y_i), i =1, \ldots, n$ from the model, the log-likelihood for $(\beta, s) \in \mathbb R^d \times (0, \infty)$ is given by
 \begin{eqnarray*}
\ell_n(\beta, s)  = \sum_{i=1}^n \log  f^X(X_i) - n \log s +  \sum_{i=1}^n \log f\left(\frac{Y_i-\beta^T X_i}{s}\right). 
\end{eqnarray*}
Since the term $\sum_{i=1}^n \log  f^X(X_i)$ does not depend on the unknown parameters, we will denote again by $\ell_n(\beta, s)$  the \lq\lq conditional \rq\rq \ log-likelihood and define the MLE of $(\beta_0, s_0)$  as
\begin{eqnarray*}
(\widehat{\beta} _n, \widehat s_n) & \coloneqq  &  \arg \max_{(\beta, s) \in \mathbb{R}^d  \times (0, \infty)} \ell_n(\beta, s) \\
& = & \arg \max_{(\beta, s) \in \mathbb{R}^d  \times (0, \infty)} \left \{ - n \log s +  \sum_{i=1}^n \log f\left(\frac{Y_i-\beta^T X_i}{s}\right)  \right \}.
\end{eqnarray*}

Before investigating convergence properties of the MLE, we start with the following theorem which shows existence of this estimator. Without further assumptions on the reference density $f$, showing existence of an MLE is in general a hard task.  The main difficulty lies in the fact that the domain of maximization is not compact. In the sequel, we will make the following assumption about the tail behavior of the density $f$:  
\begin{description}
\item \textbf{(A4)} There exist $\alpha > 0$ and $C > 0$ such that $f(t)  \le C \exp(- \vert t \vert^\alpha)$. 
\end{description}
Note that this assumption is satisfied by the Laplace ($\alpha =1$) and Gaussian $(\alpha =2)$ distributions.  In Section \ref{eff} below, we will consider the case where $f$ is exactly given by $d_\alpha \exp(- \vert t \vert^\alpha)$ for some normalizing constant $d_\alpha > 0$ that can be explicitly determined using the Gamma function.

\medskip

\begin{theorem}\label{exist}
Assume that (A4) holds. For any fixed $n \ge 1$, we assume further that
\begin{description}
\item [(i)] there exists no vector $v \in \mathbb R^d  \setminus \{0\}$ such that $v^T X_i =0$ for all $i =1, \ldots, n$,
\item [(ii)] there exists no vector $\beta \in \mathbb R^d$ such that $Y_i = \beta^T X_i$  for all $i =1, \ldots, n$
\end{description}
Then, there exists at least one MLE  $(\widehat{\beta}_n,\widehat s_n) \in \mathbb{R}^d \times (0, \infty)$.
\end{theorem}

\medskip

Note that (ii) hinders having the likelihood diverge to $\infty$. In fact, if there exists $\widetilde{\beta} \in \mathbb R^d$ such that $Y_i = \beta^T X_i, i=1, \ldots, n$ then
$$
\frac{1}{n} \ell_n(\widetilde{\beta}, s) = -\log(s) + \log(f(0))  \nearrow \infty, \ \ \textrm{as $s \searrow 0$}.
$$

\subsection{Asymptotic normality of the MLE}

\subsubsection{Consistency}
In this section, we will show that any MLE is consistent.  Recall that we could exhibit $R_0$, $a_0$ and $b_0$ (which depend on $\alpha$ and the data $(X_i, Y_i), i=1, \ldots, n$) such that the maximization of the log-likelihood is shown to be only meaningful on the compact set $\mathcal{B}(0, R_0) \times [a_0, b_0]$.  We will show the existence of a compact that does not depend on the data, and on which the maximization task can be restricted.  In the sequel, we will assume that 
\begin{description}
\item  \textbf{(A5)}  $\displaystyle \mathbb E_{(\beta_0, s_0)} \left \vert \log f \left(\frac{Y- \beta^T_0  X}{s_0} \right) \right \vert < \infty$. 

\item  \textbf{(A6)}  $\mathbb E[ \Vert X \Vert^{2\alpha \vee 2}  ]  < \infty$. 

\item  \textbf{(A7)}  The density  $f$ is monotone (non-increasing) on $[0, \infty)$.

\item  \textbf{(A8)}  $ E_{(\beta_0, s_0)} \left[ \log f \left(\frac{\vert Y \vert + R \Vert X \Vert }{a}  \right) \right]^2 < \infty $ for any $R > 0, a > 0$.

\end{description}

\bigskip

Some remarks are in order. Since $f$ is even, Assumption (A7) implies that $f$ is monotone non-decreasing on $(-\infty, 0]$. It is worth noting that for the theory to hold, this assumption can be relaxed by assuming that $f$ changes monotonicity only a finite number of times.  In fact, monotonicity is an important element in the proof of uniform consistency since it preserves the VC property. Assumption (A8) ensures that the class of functions $(x, y) \mapsto \log f( (y- \beta^T x)/s) $, as $(\beta, s) \in \mathcal{B}(0, R^*) \times [a^*, b^*]$,  admits an envelope that has a finite second moment with respect to the true distribution.  Note that this assumption is satisfied in case $f(t)  =  d_\alpha \exp(- \vert t \vert^\alpha)$ for some $\alpha > 0$. In fact, we have in this case that
\begin{eqnarray*}
\left \vert  \log f \left(\frac{\vert y \vert + R \Vert x \Vert}{a}  \right) \right \vert  & \le   &  \vert \log d_\alpha \vert + \frac{1}{a^{\alpha}} (\vert y \vert+  R \Vert x \Vert)^{\alpha} \\
& \le & \vert \log d_\alpha \vert + \frac{2^{\alpha-1} \vee 1}{ a^\alpha} (\vert y \vert^\alpha +  R^\alpha \Vert x \Vert^\alpha)
\end{eqnarray*}
and hence
\begin{eqnarray*}
&& \mathbb E_{(\beta_0, s_0)} \left[\left \vert  \log f \left(\frac{\vert Y \vert + R \Vert X \Vert}{a}  \right) \right \vert \right]^2 \\
&& \le     2 \vert \log d_\alpha \vert^2 +  4 \frac{2^{2(\alpha-1)} \vee 1}{a^{2\alpha}} 
 \big(\mathbb E_{(\beta_0, s_0)}[\vert Y \vert^{2 \alpha}] + R^{2\alpha}   \mathbb E[\Vert X \Vert^{2\alpha}]\big) < \infty
\end{eqnarray*}
by Assumption (A6) and the fact that 
$$
\mathbb E_{(\beta_0, s_0)} [\vert Y \vert^{2 \alpha}  ] \le (2^{2\alpha-1} \vee 1) \big(s^{2\alpha}_0\mathbb E(\vert \epsilon \vert^{2\alpha}) + \Vert \beta_0 \Vert^{2\alpha} \mathbb E(\Vert X \Vert^{2\alpha} )\big).
$$
\bigskip

\par \noindent We start with the following proposition.

\begin{proposition}\label{Determ}
There exist deterministic constants $R^* > 0, a^* > 0$ and  $b^* > 0$  such that with probability 1 the problem of maximizing the log-likelihood can be restricted to the compact set  $\mathcal{B}(0, R^*) \times [a^*, b^*]$ for $n$ large enough.
\end{proposition}

Now, we are ready to state the consistency result. As we first need to establish uniform consistency of the log-likelihood on the relevant space of parameters, the proof of the next theorem relies on arguments from empirical process theory.  For more details, see Appendix B.

\begin{theorem}\label{Consis}
 Let $(\widehat \beta_n, \widehat s_n)$ denote a maximum likelihood estimator (MLE). Then, under the assumptions above
$$
(\widehat \beta_n, \widehat s_n) \to_{\mathbb P} (\beta_0, s_0).
$$
 \end{theorem}

\bigskip

\subsubsection{Weak convergence}

In this section, we will use the consistency established above to show that any MLE $(\widehat \beta_n, \widehat s_n)$ is asymptotically normal and exhibits the asymptotic variance. To derive this result, we need the following assumptions. Below, $f_s$ denotes $1/s f(\cdot/s)$.
\begin{description}
\item \textbf{(A9)}  The density $f$ is twice continuously differentiable on the interior of its support. 


\item \textbf{(A10)}  For all $\delta > 0$, there exists $\eta > 0$ such that if $\Vert \beta^* - \beta_0 \Vert + \vert s^* - s_0 \vert \le \eta $ it holds for all $1 \le i, j \le d$
\begin{eqnarray*}
\sup_{(x, y)} \left \vert  \frac{\partial^2 \log f_s(y-\beta^Tx)}{\partial \beta_i \beta_j }|_{(\beta, s) = (\beta_0, s_0)} - \frac{\partial^2 \log f_s(y-\beta^Tx)}{\partial \beta_i \beta_j }|_{(\beta, s) = (\beta^*, s^*)}  \right \vert  < \delta
\end{eqnarray*}

\item \textbf{(A11)}  The matrix $\mathbb E[X X^T]$ is positive definite, and 
$$
\int t^2 \frac{(f'(t))^2}{f(t)} dt  \in (0, \infty).
$$
\end{description}

\medskip

\begin{theorem}\label{asympMLE}
Under the assumptions above, the MLE $\widehat \theta_n = (\widehat \beta_n, \widehat s_n)$ is asymptotically normal. More specifically, if $\theta_0 = (\beta_0, s_0)$ then
\begin{eqnarray*}
\sqrt n (\widehat \theta_n - \theta_0)  \to_d \mathcal N(\mathds{O}_{d+1}, I^{-1}_0)
\end{eqnarray*}
where
$$I_0 =\left( \begin{array}{cc}
       c_1\mathbb{E}[XX^T] & \mathds{O}_d \\
        \mathds{O}^T_d & c_2
    \end{array}
\right) $$
with
\begin{eqnarray*}
    \mathds{O}_k &\coloneqq& (0,0,...,0)^T \in \mathbb{R}^k\\
    c_1 &\coloneqq& \frac{1}{s_0^2}\int_{\mathbb{R}}\frac{(f'(t))^2}{f\left(t\right)}dt\\
    c_2 &\coloneqq& \frac{1}{s_0^2}\left(\int_{\mathbb{R}}\frac{(f'(t))^2}{f\left(t\right)}t^2 dt -1 \right).
\end{eqnarray*}

\end{theorem}

\bigskip

\begin{remark}
The constant $c_2$ is always strictly positive. In fact, using the Cauchy-Schwarz inequality, we can write
\begin{eqnarray*}
\left(\int t f'(t) dt\right)^2 & = &  \left(\int t \frac{f'(t)}{f(t))} f(t) dt   \right)^2 \\
& \le & \left(\int t^2 \frac{(f'(t))^2}{f(t))} dt   \right) \int f(t) dt = \int t^2 \frac{(f'(t))^2}{f(t)} dt.
\end{eqnarray*}
Since $\int t f'(t) dt = -1$ (see the proof of Proposition \ref{asybehavols}), it follows that $c_2 \ge 0$. By the characterization of the equality part in the  Cauchy-Schwarz inequality, $c_2 = 0$ if and only if there exists $\lambda \in \mathbb R$ such that for a.e. $t$
\begin{eqnarray*}
t \frac{f'(t)}{f(t))}  = \lambda
\end{eqnarray*}
or equivalently
\begin{eqnarray*}
\log(f(t))  =  \lambda \log(\vert t \vert) + c, \ \text{for some $c \in \mathbb R$}.
\end{eqnarray*}
This means that $f(t)  = \exp(c) \vert t \vert^{\lambda}$. Since $t \to \vert t \vert^{\lambda}$ is not integrable over $\mathbb R$, we conclude that $c_2 > 0$.

\end{remark}

\bigskip


\section{Asymptotic relative efficiency with respect to the OLSE}\label{sec: AsympEffic}

\subsection{Asymptotic efficiency}
In this section, we prove that the MLE for $\beta_0$, derived from the joint MLE of the pair $(\beta_0,s_0)$, is asymptotically at least as efficient as the OLSE, denoted in the sequel by $\widehat \beta^{\text{OLS}}_n$.  Under Assumption (A12), we can show that 
\begin{eqnarray}\label{asympOLS}
\sqrt{n}(\widehat{\beta}^{\text{OLS}}_n-\beta_0) \xrightarrow[]{d} \mathcal{N}_d\left(0,\sigma^2 (\mathbb{E}[XX^T])^{-1}\right). 
\end{eqnarray}
A proof of the convergence in (\ref{asympOLS}) is provided in Appendix B for the sake of completeness. Using (\ref{asympOLS}) we can now show the following proposition.

\medskip

\begin{proposition}\label{asybehavols}
 The MLE of $\beta_0$, $\widehat \beta_n$, is at least as efficient as the OLSE. More specifically, if $\eta$ denotes the asymptotic relative efficiency of $\widehat \beta_n$ with respect to $\widehat{\beta}^{\text{OLS}}_n$, then 
 \begin{eqnarray}\label{eta}
 \eta =  \left[\left(\int t^2 f(t) dt\right) \left(\int \frac{(f'(t))^2}{f(t)} dt \right)\right]^{-1} \le 1
 \end{eqnarray}
 with equality if and only if $\epsilon \sim \mathcal{N}(0, \sigma^2)$. 
\end{proposition}

\subsection{Quantifying the gain for using the MLE under certain noise distributions}\label{eff}

In this section, we consider three families of symmetric densities for the noise $\epsilon$. They are parameterized by $\gamma > 0$ in different ways with the common point that a density $f_{\gamma, i}$ in each of the three families satisfies $\int_{\mathbb R} t^2 f_{\gamma, i}(t) dt = 1$ for $i=1, 2, 3$.

\begin{enumerate}

\item $f_{\gamma,1}(t) = \displaystyle d_\gamma \exp\left( - c_\gamma \vert  t \vert^{\gamma}  \right)$. In this case 
\begin{eqnarray*}
c_\gamma  =  \left(\frac{\Gamma\left(\frac{3}{\gamma}\right)}{\Gamma\left(\frac{1}{\gamma}\right)}\right)^{\gamma/2}, \ \text{and} \ d_\gamma =   \frac{\gamma}{2} \left(\frac{\Gamma\left(\frac{3}{\gamma}\right)}{\Gamma\left(\frac{1}{\gamma}\right)^{3}}\right)^{1/2}.
\end{eqnarray*}

\vspace{0.1cm}
\item  $f_{\gamma,2}(t)   =  \displaystyle d_\gamma  \vert t \vert^{\gamma -1} \exp(- c_\gamma \vert t \vert), \gamma \ge 1$.  In this case, 
\begin{eqnarray*}
c_\gamma =  \left(\frac{\Gamma(\gamma +2)}{\Gamma(\gamma)}\right)^{1/2}, \ \ \text{and} \ \  d_\gamma =   \frac{1}{2 \Gamma(\gamma)} \left(\frac{\Gamma(\gamma+2)}{\Gamma(\gamma)}\right)^{\gamma/2}.  
\end{eqnarray*}

\vspace{0.1cm}

\item $f_{\gamma,3}(t)  = \displaystyle   d_\gamma \vert t \vert^{\gamma -1}  \exp(- c_\gamma \vert t \vert^\gamma), \gamma \ge 1$.  In this case, 
\begin{eqnarray*}
c_\gamma  =  \Gamma\left(\frac{2}{\gamma} +1 \right)^{\gamma/2}, \ \ \text{and} \ \  d_\gamma =   \frac{\gamma}{2}  \Gamma\left(\frac{2}{\gamma} +1 \right)^{\gamma/2}.  
\end{eqnarray*}
\end{enumerate}
Above,  $\Gamma(z)  = \int_0^\infty t^{z-1} \exp(-t) dt$ for $z\in (0, \infty)$,

\medskip
Note that Laplace-families have $\gamma=1$ for families 1-3 and Gaussian have $\gamma=2$ for the first family.  We would like to find an optimal $\gamma^*$ for which the relative efficiency between OLSE and MLE is maximized.

For $i=1, 2, 3$, let $\gamma \mapsto \eta_i(\gamma)$ denote  the asymptotic relative efficiency function of the MLE with respect to the OLSE for each of the three families. To find the expressions of $\eta_i(\gamma)$, we need to compute the integral $\int (f'_{\gamma, i}(t))^2/ f_{\gamma, i}(t) dt )dt$ since $\int t^2 f_{\gamma, i}(t) dt = 1$. When the  integration is possible, one has to resort to some involved algebra, use of the change of variable $t \mapsto t^\gamma = x$ (for the first and third families) and the well-known property of the Gamma function: $\Gamma(a) = (a-1) \Gamma(a-1)$ for $a > 1$. We find that
\begin{eqnarray*}
\eta_1(\gamma) =  \frac{\Gamma\left(\frac{1}{\gamma}\right)^2}{\gamma^2 \ \Gamma\left(\frac{3}{\gamma}\right) \Gamma\left(2-\frac{1}{\gamma}\right)}, \  \  \gamma \in (1/2, \infty),
\end{eqnarray*}
\begin{eqnarray*}
\eta_2(\gamma) =  \frac{\Gamma(\gamma) (\gamma -2)}{\Gamma(\gamma +2)}, \  \  \gamma \in (2, \infty),
\end{eqnarray*}
and 
\begin{eqnarray*}
\eta_3(\gamma) =  \frac{1}{\Gamma\left(\frac{2}{\gamma} +1 \right) \Gamma\left(1-\frac{2}{\gamma}\right) (\gamma -1)^2}, \  \  \gamma \in (2, \infty).
\end{eqnarray*}

\begin{figure}[H]
    \centering
\includegraphics[width=0.7\textwidth]{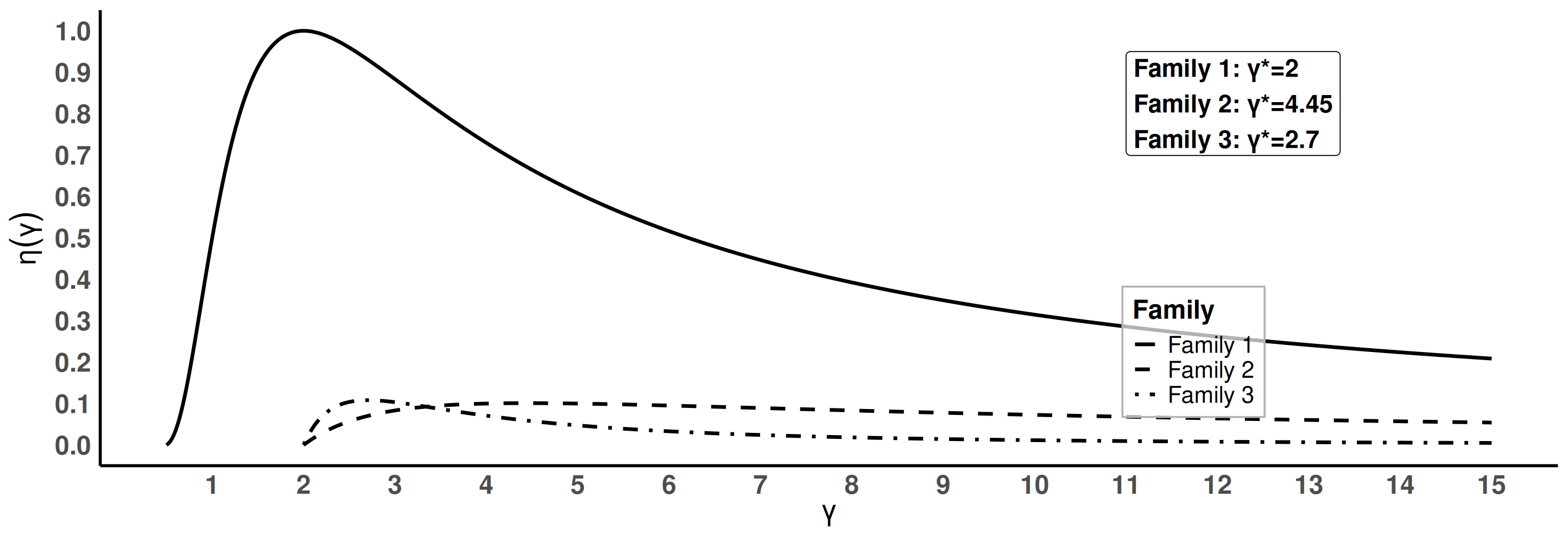}    
    \caption{Asymptotic relative efficiency of the MLE versus the OLSE for the three families as a function of $\gamma$. $\gamma^*$ are the values for which it is maximized.}
    \label{fig:eta}
\end{figure}

\subsection{When the noise is asymmetric}
In most regression settings, the noise is assumed to be symmetric around $0$. This is also the working assumption under which the statistical gain of the MLE over the OLS is derived.  One referee raised the question of whether the gain would be affected if the noise distribution is asymmetric.  First of all, note that when the regularity assumptions under which asymptotic normality of the MLE is valid, and the noise distribution is well-specified, the MLE should be the most efficient estimator. Hence, the MLE is expected to be at least efficient as any other estimator.  However, deriving the theoretical gain seems to be a much harder endeavor in this case. In particular, the ratio of the determinants of the covariance matrices of the MLE and the OLS is not an appropriate measure in this case.  In fact, let $c_1$ and $c_2$ be the same constants as above in Theorem \ref{asympMLE}. Additionally, consider the constants
\begin{eqnarray*}
c_3  =  \frac{1}{s_0^2} \int t \frac{(f'(t))^2}{f(t)} dt \end{eqnarray*}
and  
$$\kappa = \frac{c_3^2}{c_1 c_2}$$
in case $c_2 \ne 0$. Finally, let $\mu = \mathbb E(X)$ and $V =  \mathbb E[X X\top] - \mu \mu^\top$ the covariance matrix of $X$, assumed to be invertible.  \\

After some tedious algebra, we can show that the ratio of the determinants raised to the power $1/d$ is given by  \begin{eqnarray}\label{Gainnonsym}
\eta =  \frac{1}{c_1 \sigma^2} \left[\frac{1 + \mu^\top V^{-1} \mu}{1 + (1-\kappa) \mu^\top V^{-1} \mu}  \right]^{1/d}.
\end{eqnarray}
Note that in the symmetric case, $c_3 = \kappa = 0$, and $\eta$ is always $\le 1$. In the asymmetric case, we obtain the same gain as above if $\mu = 0$. However, the value of $\eta$ in (\ref{Gainnonsym}) is not necessarily below $1$ and can even take values much larger than 1 depending on $\mu^\top V^{-1} \mu$. This does not contradict the theory but rather reveals that the ratio of the determinants is not always a good measure for comparing the efficiency of two estimators.  In fact, the determinant of a covariance matrix measures the volume of the confidence ellipsoid and is multiplicative.  When the noise density is asymmetric, we know that the matrix of the asymptotic covariance of the MLE is no longer proportional to that of the OLSE. The ratio of the determinants might be heavily affected by the directions where the MLE is doing slightly worse than the OLSE.

In the next section, we investigate the statistical gain in case the noise distribution belongs to the scaled family of shifted exponential distribution using other comparison measures which remain more faithful to the overall (better) performance of the MLE even if it is a bit worse than the OLSE along some directions.

\section{Simulations and a data application}\label{simu}
In this section, all figures were made with the statistical software \texttt{R}.
\subsection{Synthetic data}
We now study the performance of the MLE and the OLSE through synthetic regression experiments. The focus is on regression models with error terms from the three symmetric distribution families introduced in Section~\ref{eff}. Across $M=20000$ replications and for varying sample sizes $n$, dimensions $d$, and shape parameters $\gamma$, we proceed in three steps. First, we illustrate how the bias and variance of the two estimators differ. Second, we further look at the effect of the dimension $d$ on the relative performance of the estimators as measured by their mean squared distance to $\beta_0$. Finally, we examine simulated and theoretical values of the asymptotic relative efficiency $\eta$ defined in equation \eqref{eta}.
The code use to obtain the results shown in the coming section can be found under GitHub repository \url{https://github.com/juslecl/PubThesis/tree/master}.
\subsubsection{Data generation and estimation procedure}

To simulate data under each of the three distributions introduced in section~\ref{eff}, we consider sampling procedures that involve existing standard distributions implemented in R, and which we further rescale to exactly match our target densities. The transformation parameters are derived analytically to ensure that the generated values follow the desired distribution precisely.

\paragraph{Generalized Normal Distribution (\(f_1^\epsilon\))} Let \( X \sim \text{GN}(0,1,\gamma) \) be generated via \texttt{gnorm::rgnorm()}, with density
\[
f_{GN}(x) = \frac{\gamma}{2\Gamma(1/\gamma)} \exp\left( -|x|^\gamma \right),
\]
We define
\[
c_\gamma := \left( \frac{\Gamma(3/\gamma)}{\Gamma(1/\gamma)} \right)^{\gamma/2},
\]
and define \( T = s \cdot c_\gamma^{-1/\gamma} X \). The resulting variable has density
\[
f_T(t) = \frac{1}{s} \cdot d_\gamma \exp\left( -c_\gamma \left| \frac{t}{s} \right|^\gamma \right), \quad
d_\gamma = \frac{\gamma}{2} \left( \frac{\Gamma(3/\gamma)}{\Gamma(1/\gamma)^3} \right)^{1/2}.
\]

\paragraph{Generalized Gamma with Linear Exponent (\(f_2^\epsilon\))}

Let \( X \sim \text{Gamma}(\gamma, c_\gamma^{-1}) \) be generated via \texttt{stats::rgamma()},
\[
f_X(x) = \frac{c_\gamma^\gamma}{\Gamma(\gamma)} x^{\gamma-1} \exp(-c_\gamma x), \quad x > 0,
\]
with constants,
\begin{eqnarray*}
c_\gamma =  \left(\frac{\Gamma(\gamma +2)}{\Gamma(\gamma)}\right)^{1/2}, \ \ \text{and} \ \  d_\gamma =   \frac{1}{2 \Gamma(\gamma)} \left(\frac{\Gamma(\gamma+2)}{\Gamma(\gamma)}\right)^{\gamma/2}.  
\end{eqnarray*}
To symmetrize, we draw \( T = \pm sX \) with equal probability. The resulting density is
\[
f_T(t) = \frac{1}{s} \cdot d_\gamma \cdot |t|^{\gamma - 1} \exp\left( -c_\gamma \cdot \frac{|t|}{s} \right).
\]

\paragraph{Generalized Weibull Distribution (\(f_3^\epsilon\))}

Let \( X \sim \text{Weibull}(\gamma, c_\gamma^{-1/\gamma}) \) be generated via \texttt{rweibull()}, with density
\[
f_X(x) = \frac{\gamma}{c_\gamma^{-1/\gamma}} \left( \frac{x}{c_\gamma^{-1/\gamma}} \right)^{\gamma - 1} \exp\left( -\left( \frac{x}{c_\gamma^{-1/\gamma}} \right)^\gamma \right),
\]
and constants
\[
c_\gamma = \Gamma\left( \frac{2}{\gamma} + 1 \right)^{\gamma/2}, \quad
d_\gamma = \frac{\gamma}{2} \Gamma\left( \frac{2}{\gamma} + 1 \right)^{\gamma/2}.
\]
We again symmetrize by setting \( T := \pm sX \) with equal probability, yielding the density
\[
f_T(t) = \frac{1}{s} \cdot d_\gamma \cdot |t|^{\gamma - 1} \exp\left( -c_\gamma \cdot \frac{|t|^\gamma}{s} \right).
\]

\paragraph{Simulation setup.} \label{simsetup} 
For each parameter configuration $(n,d,\gamma)$ and each family $f^\epsilon_j$, $j\in\{1,2,3\}$, we generate $M=20000$ independent datasets. The design matrix $X\in\mathbb{R}^{n\times d}$ has i.i.d.\ rows from $\mathcal{N}_d(\mu_X,\sigma_X^2 I_d)$. Parameters are drawn once with seed number 2105 and fixed across replications:
\[
\beta_0 \sim \mathcal{U}_d[-5,5],\quad 
\mu_X \sim \mathcal{U}_d[-3,3]
\]
The standard deviation of $X$ and $\epsilon$ are fixed to 
\[ \sigma_X = 2,\quad 
s_0 =2.75
\]
 respectively. 

\subsubsection{Computation and performance evaluation}
Given data \( \{(x_i, y_i)\}_{i=1}^n \), we compute the maximum likelihood estimates (MLE) of the regression parameter \( \beta_0 \in \mathbb{R}^d \) and the scale parameter \( s_0 > 0 \) by minimizing the negative log-likelihood of the assumed error distribution. Optimization uses the Broyden-Fletcher-Goldfarb-Shanno (BFGS) algorithm (\texttt{optimx::optimx()}), initialized with the OLSE for $\beta_0$ and the empirical standard deviation of OLSE residuals for $s_0$. We also compared a maximization procedure where different initializing values for the standard deviation are picked from a grid consisting of values strictly smaller, equal and larger than the empirical standard deviation of OLSE residuals. These additional experiments were motivated by a referee's remark about the fact that initializing the optimization procedure with the empirical standard deviation of the OLSE residuals may not be appropriate for heavy-tailed or skewed distributions such as family \#2 and \#3.  We found that  the resulting parameter estimates differed only at the order of $10^{-4}$, and no substantial changes in the likelihood values were observed. Given these findings, together with the simplicity and ease of implementation of our approach, we chose to retain it in the present work.

We would like to note the connection to GAMLSS framework, to which another reviewer drew our attention; see \cite{rigby2005}. Our regression model can indeed be fitted in the GAMLSS framework if we choose the location to be $\mu(x) = \beta_0\top x$, a constant scale and fixed shape. Out of curiosity, we attempted to recompute the MLEs from our simulations using the R package \textbf{gamlss}. This package provides a wide range of pre-implemented standard distributions and allows users to define custom \texttt{gamlss.family} objects. Since the distributions corresponding to families \#2 and 3 considered in the present work were not available, we tried to implement them manually within the \texttt{gamlss} framework. However, after many attempts, we were unable to successfully define the corresponding \texttt{gamlss.family} objects, as systematic implementation errors were returned.
One practical advantage of our proposed approach is that it is  considerably easy to implement using existing R packages and requires only limited additional coding expertise.  Also, note that in this paper we  explicitly provide asymptotic comparisons of the MLE versus the OLSE under the “known up to scale” noise distribution, and quantify the “reward” as a function of deviation from Gaussianity.  GAMLSS is primarily a general modeling and fitting framework, and in this sense it does not by itself supply the same theoretical efficiency analysis specialized to our setting. 

\subsubsection{Bias and variance}

We consider diverse scenarios and display boxplots of the MLE's and OLSE's mean Euclidean distance to the true $\beta_0$ over $M=20000$ replications.

The MLE uniformly outperforms the OLSE in terms of both bias (median error) and variance (spread of the boxplots). This observation aligns with our theoretical expectations: both estimators are known to coincide in performance only when the error distribution is Gaussian. In the opposite situation, the performance gap increases significantly under heavier-tailed error distributions, i.e., from \( f^\epsilon_1 \) to \( f^\epsilon_3 \), and also as the shape parameter \( \gamma \) increases.

Finally, note that the performance of both estimators worsens as $d$ increases, and we do not observe a widening performance gap between MLE and OLSE as the dimension gets larger.

\begin{figure}[H]
    \centering
    \includegraphics[width=\linewidth]{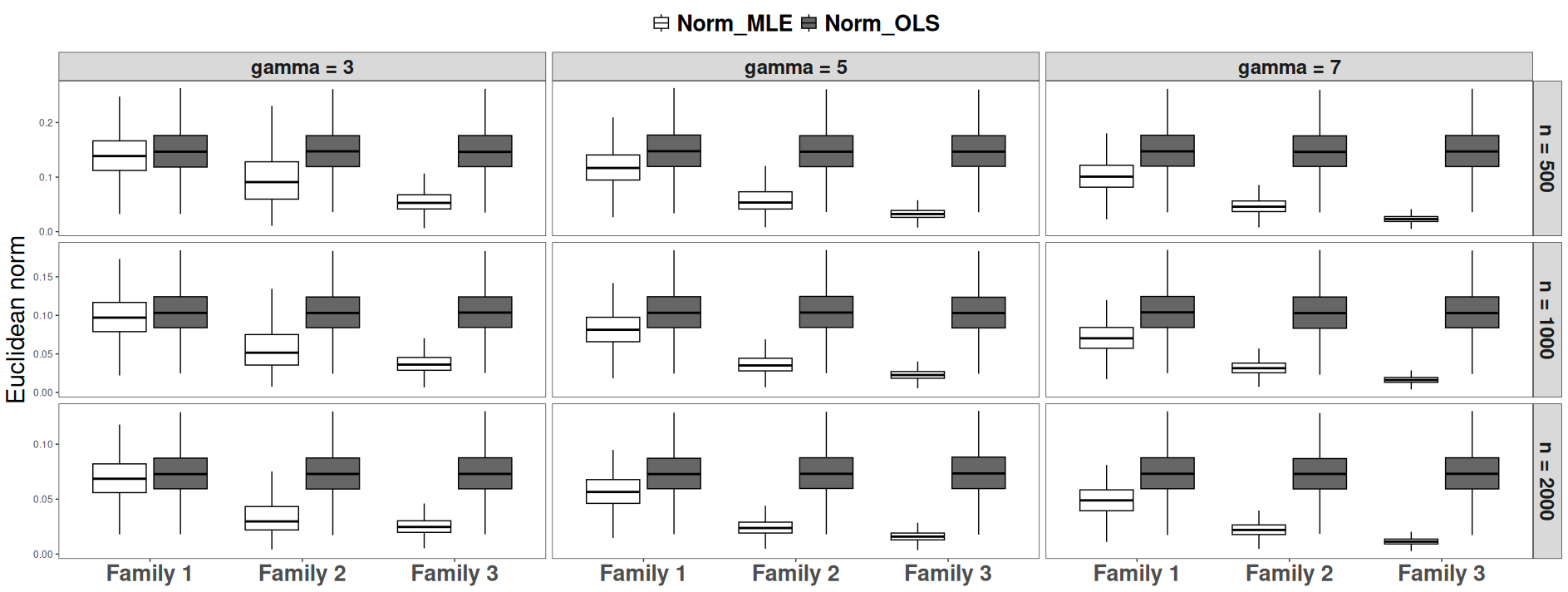}
    \caption{Boxplot for the Euclidean norm from MLE resp. OLS to true $\beta_0$, for $d=3$, $n\in \{500,1000,2000\}$ $\gamma \in \{3,5,7\}$}
    \label{fig:placeholder}
\end{figure}

\begin{figure}[H]
    \centering
    \includegraphics[width=\linewidth]{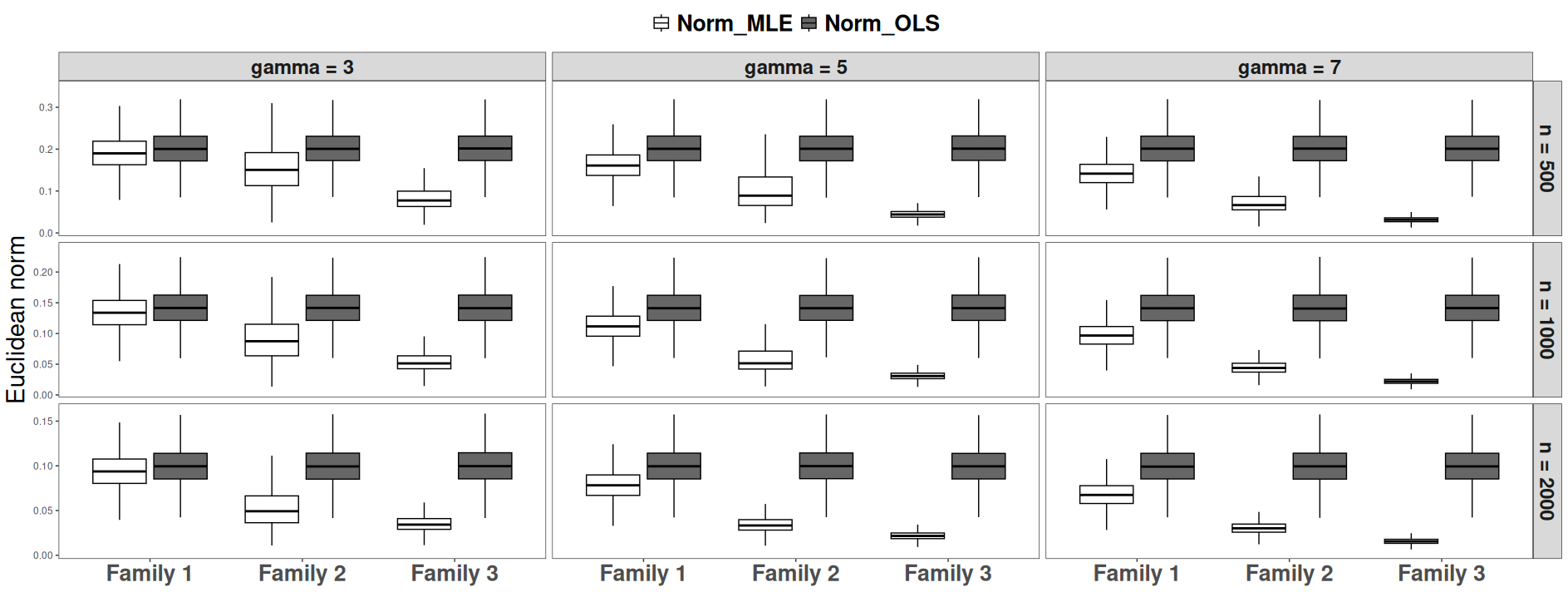}
    \caption{Boxplot for the Euclidean norm from MLE resp. OLS to true $\beta_0$, for $d=13$, $n\in \{500,1000,2000\}$ $\gamma \in \{3,5,7\}$}
    \label{fig:placeholder}
\end{figure}

\subsubsection{Asymptotic Relative Efficiency}
Finally, we compare the simulated values of $\eta$ (defined in Equation \eqref{eta}) with their theoretical values. For each scaled family, we perform 200 independent simulations, each consisting of 100 replicated datasets generated with fixed sample size $n=2000$, fixed $d=3$, and $\gamma \in \{3,5,7\}$.
For each simulation batch, we compute the estimator
\[
\hat{\eta} = \left(\frac{\det(\widehat{\text{Var}}
(\hat{\beta}_{\text{MLE}}))}{\det(\widehat{\text{Var}}
(\hat{\beta}_{\text{OLSE}}))}\right)^{1/d}
\]
where the covariance matrices are estimated empirically from the 200 replications. We then compare the resulting estimates, together with their confidence intervals, to the corresponding theoretical value of $\eta$.
\begin{figure}[H]
    \centering
    \includegraphics[width=\linewidth]{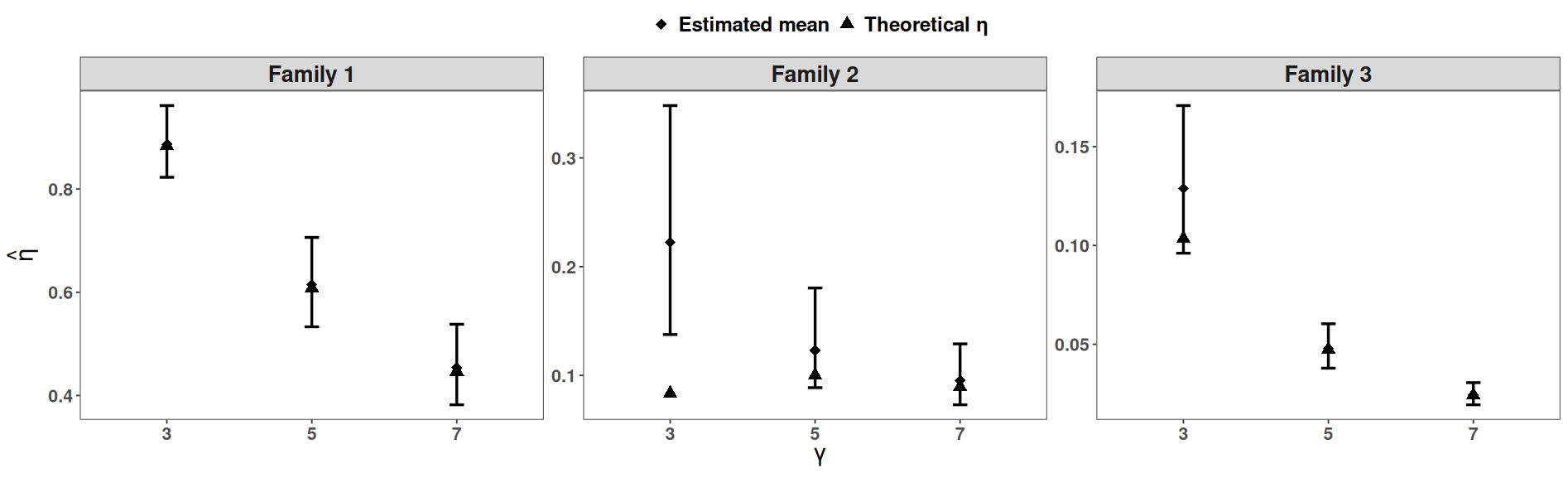}
    \caption{
Estimated mean ARE $\hat{\eta}$ vs. $\gamma$ for the three families, with 95\% percentile intervals, and theoretical $\eta$
}
    \label{are:boxplots}
\end{figure}

In Figure \ref{are:boxplots}, the estimated ARE values closely agree with the theoretical ones, supporting our derivations. An exception is observed for the second family at $\gamma=3$, where the theoretical value lies outside the 95\% percentile interval. Overall, the results confirm that $\eta$ decreases as the tails become heavier (increasing $\gamma$), with the strongest effect occurring for the third error family. \\

\medskip
The simulation results are consistent with the theoretical properties established in Section~\ref{eff}. In all non-Gaussian settings, the MLE outperforms the OLSE, with its relative advantage becoming more pronounced as the error distributions deviate further from Gaussianity, either through heavier tails or higher peakedness. Furthermore, the gap in mean squared distance widens as the dimension $d$ increases, always in favor of the MLE. This does not contradict Figure \ref{are:boxplots} where $\eta$ remains unaffected by $d$.The reason is that $\eta$ measures relative efficiency—the ratio of covariance matrix determinant—while the mean squared distance sums absolute errors across dimensions. Finally, the strong alignment between simulated and theoretical ARE values provides direct empirical validation of our efficiency analysis and further supports the use of the MLE in practice whenever the error distribution departs from Gaussianity.

\subsection{Asymmetric noise}
In the article, we suppose the noise distribution is symmetric. To investigate the impact of asymmetry on the relative performance of the MLE compared to OLS, we consider the case of the asymmetric scaled family of shifted exponential distributions: 
\begin{eqnarray*}
\left\{ f_s(x) = \frac{1}{s} f_0\!\left(\frac{x+s}{s}\right), \quad x \in \mathbb R,  \quad s > 0 \right\},
\end{eqnarray*}
where \(f_0\) is the density of an exponential distribution with intensity equal to 1. For a given scale $s > 0$, $\displaystyle f_s(x)  =  \frac{1}{s} \exp\left( - \frac{x}{s} - 1\right) \mathds{1}_{x \ge -s}$.  

The simulation setup described in \ref{simsetup} is reproduced here.

\begin{figure}
    \centering
    \includegraphics[width=\linewidth]{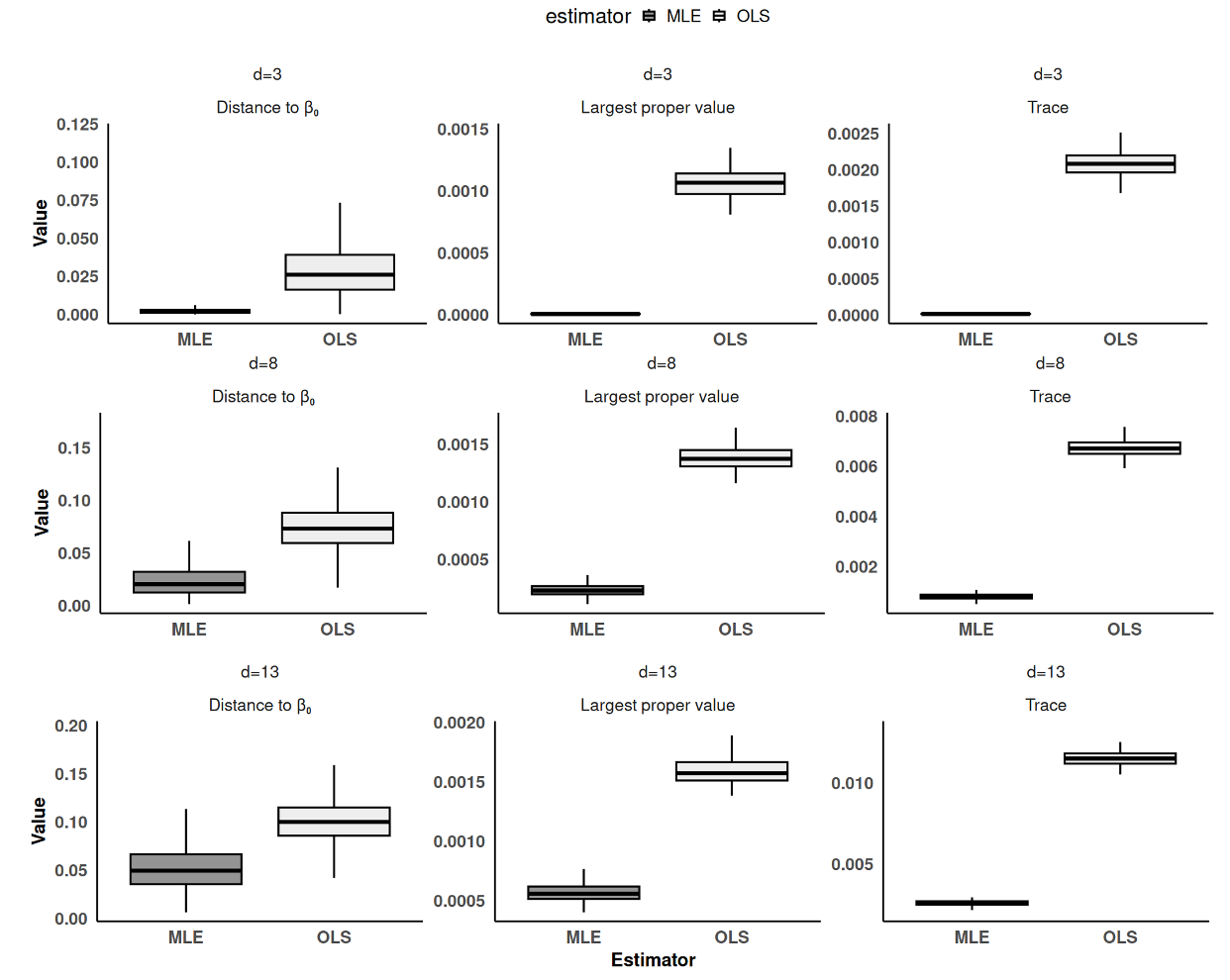}
    \caption{Boxplot for the estimator's distance to true $\beta_0$, covariance matrix largest eigenvalue, trace and for $n=2000$, $d\in \{3,8,13\}$.}
    \label{comp_exp}
\end{figure}

We compare MLE and OLS using three metrics: the trace of the covariance matrices, the largest eigenvalues of the covariance matrix, and the Euclidean distances to the true $\beta_0$. Results are plotted in Figure \ref{comp_exp}.

\begin{figure}
    \centering
    \includegraphics[width=0.7\linewidth]{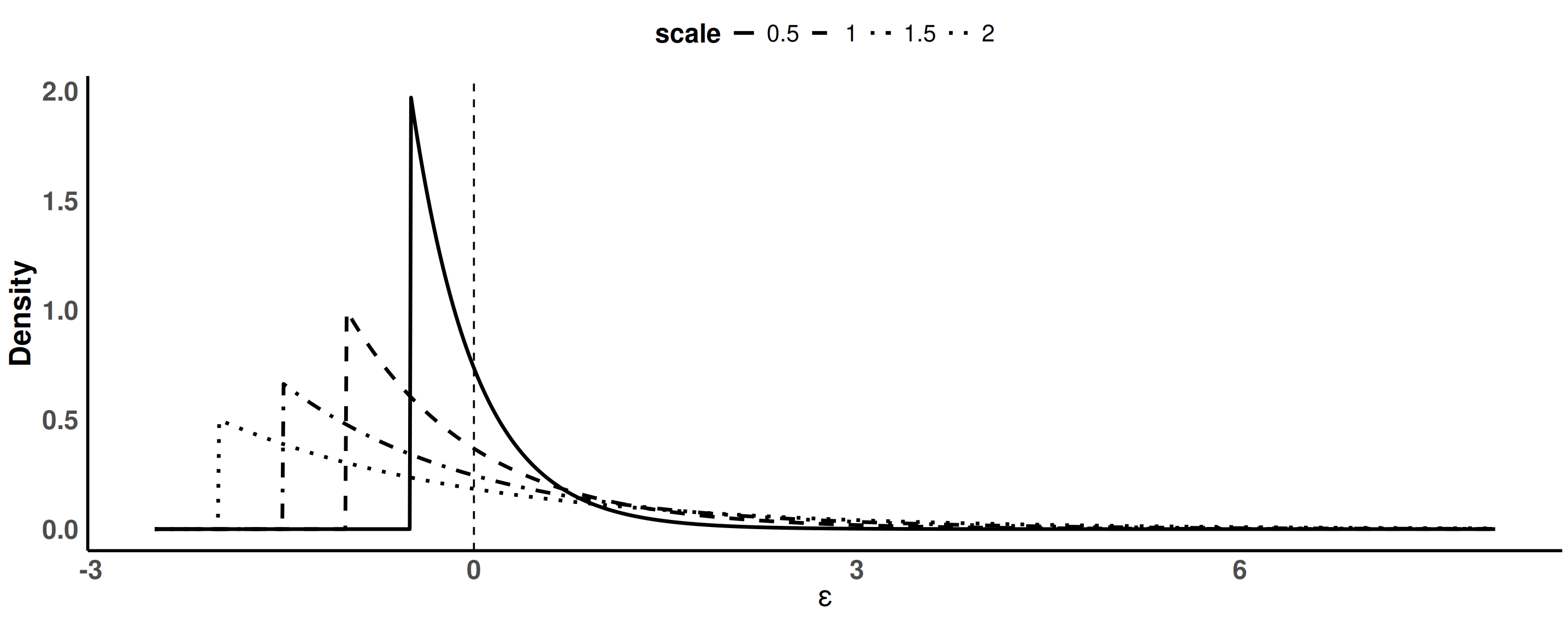}
    \caption{Densities for different values of the scale parameter \(s\).}
\end{figure}

We observe that the MLE still outperforms OLS. The substantial difference in the distributions of the largest eigenvalue and the trace of the covariance matrix, both favoring the MLE, indicates that the MLE is considerably more stable than OLS. The advantage of MLE over OLS is thus further confirmed in this asymmetric setting, and the performance gap continues to widen as the number of observations available for estimating the coefficients increases.

\subsection{Data application: The fish market dataset}

The \lq  Fish Market\rq \ dataset can be downloaded from the open website  \url{https://www.geeksforgeeks.org/machine-learning/dataset-for-linear-regression/}.  According to the description provided in the same source, the data are collected on the common fish species in fish market sales. However, the geographic location where the sales have occurred is not provided. The  dataset consists of the name of the fish species, weight (in grams), 3 length variables, height and width. Although not explicitly mentioned on the data source, the 3 length measurements should correspond to the total, fork and standard length.  There are 159 data points for 7 different Species occurring at the frequencies shown in table \ref{species}. 
\begin{table}[!h]
\begin{center}
\caption{Frequency of the fish species in the Fish Market dataset.}
\label{species}
\begin{tabular}{|cccccccc|} 
 \hline
 Species & Bream & Parkki & Perch & Pike & Roach   &  Smelt & Whitefish  \\ [0.5ex] 
 Frequency & 35  &   11   &  56  & 17  &  20 & 14 & 6
  \\ 
 \hline
 \end{tabular}
\end{center}
\end{table}

\begin{figure}[H]
\centering
    \includegraphics[width=0.6\textwidth]{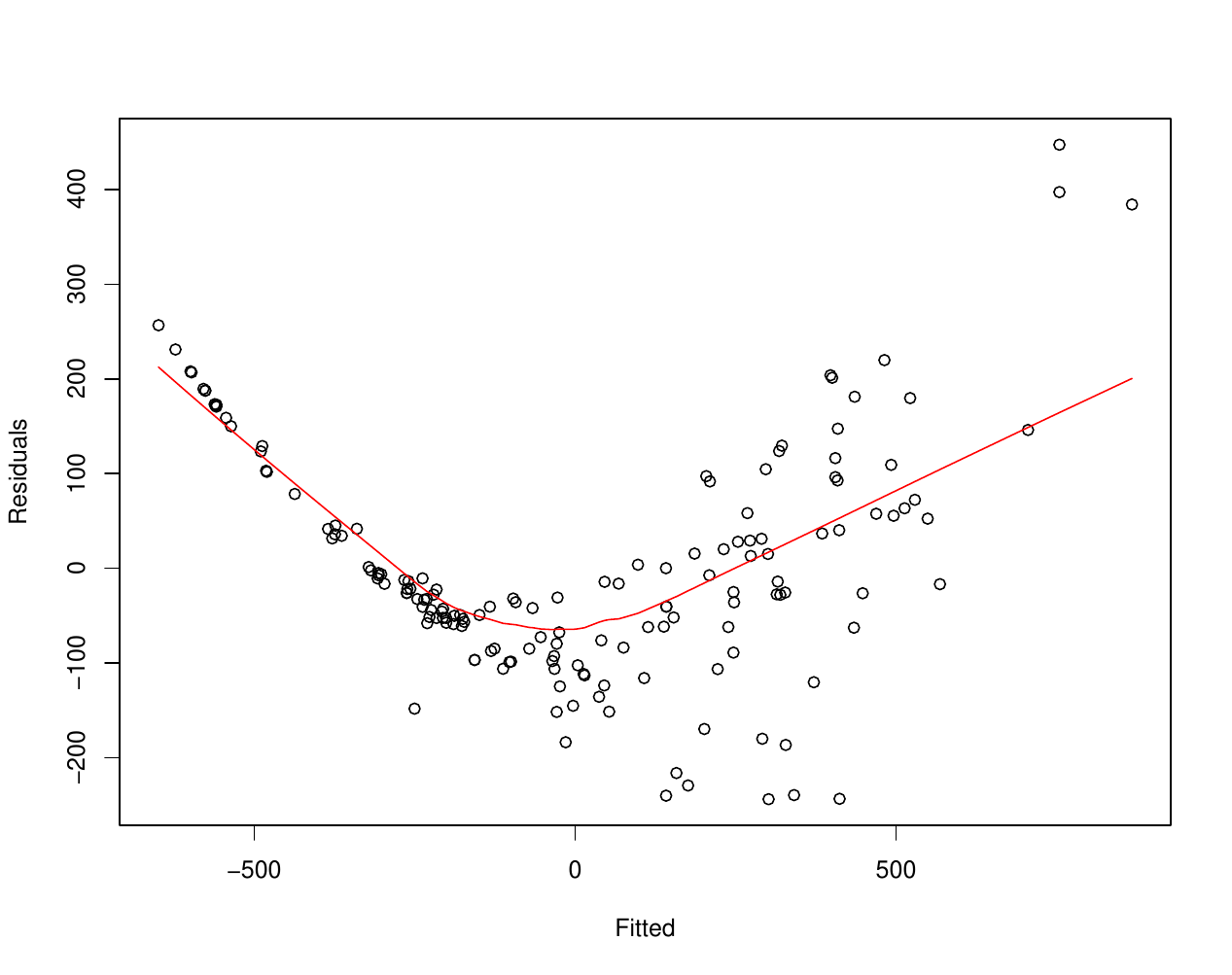}
        \caption{The residuals versus fitted values using least squares estimation in the linear model with centered fish weights as response and centered heights, widths and the 3 lengths as covariates. The added curve in red is a smoothing locally weighted polynomials obtained using the R function \text{lowess}. }
        \label{fittedresraw}
   \end{figure}

Using a multivariate linear regression, the goal is to predict fish weight based on the other variables. To conform with the assumptions under which we derived the theoretical results above, the categorical variable \lq\lq Species\rq\rq \ was not considered among the predictors. Furthermore, the response and the other (continuous) variables need to be centered since  centering allows us to consider a linear model without intercept. Figure \ref{fittedresraw} shows the residuals versus the fitted responses using a linear regression model (with all the 5 continuous and centered covariates) and the OLSE. The smoothing curve in red highlights the fact that the assumption of independence of the noise and covariates does not seem to hold.

\begin{figure}[H]
\centering
    \includegraphics[width=0.48\textwidth]{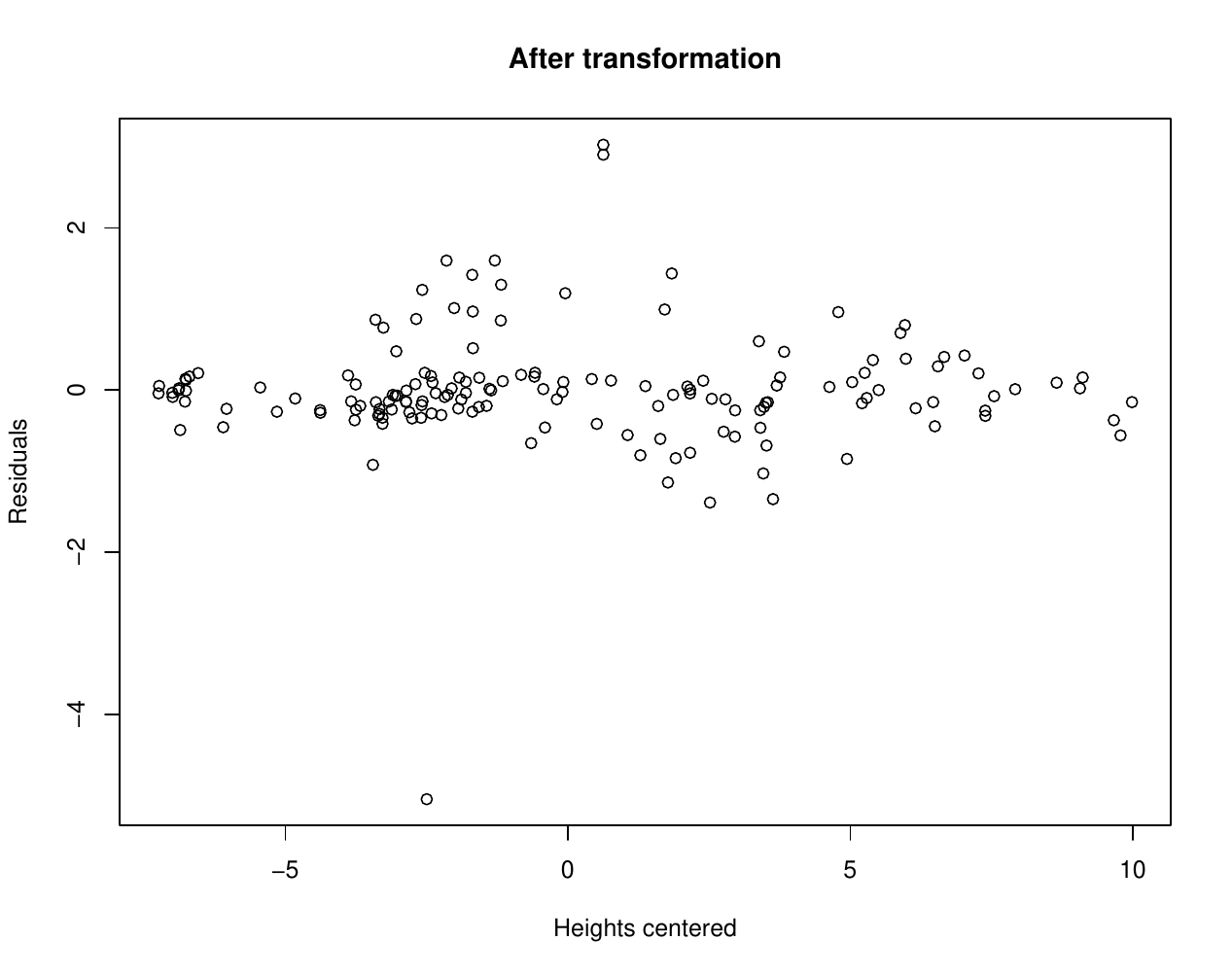}
     \includegraphics[width=0.48\textwidth]{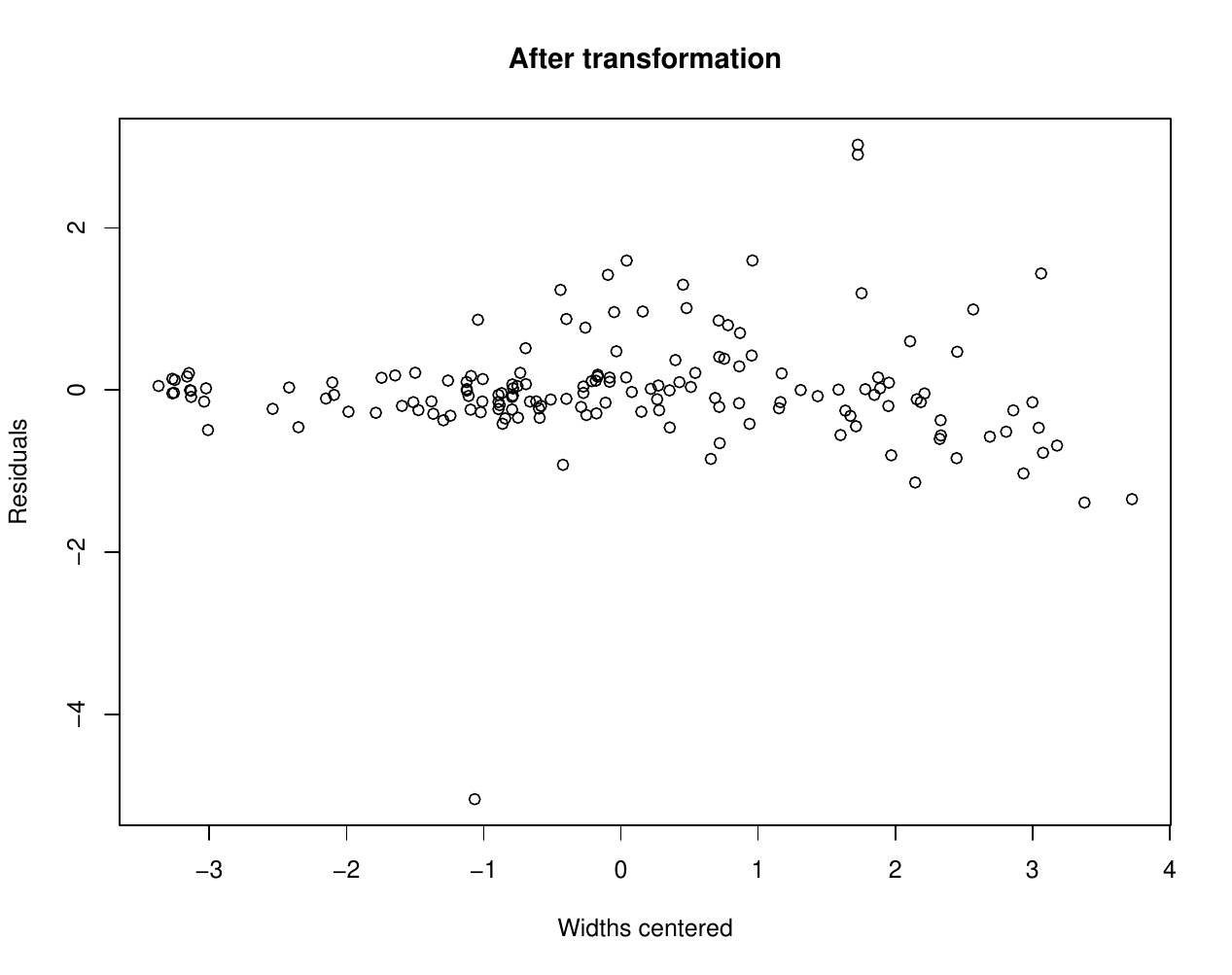}
      \includegraphics[width=0.48\textwidth]{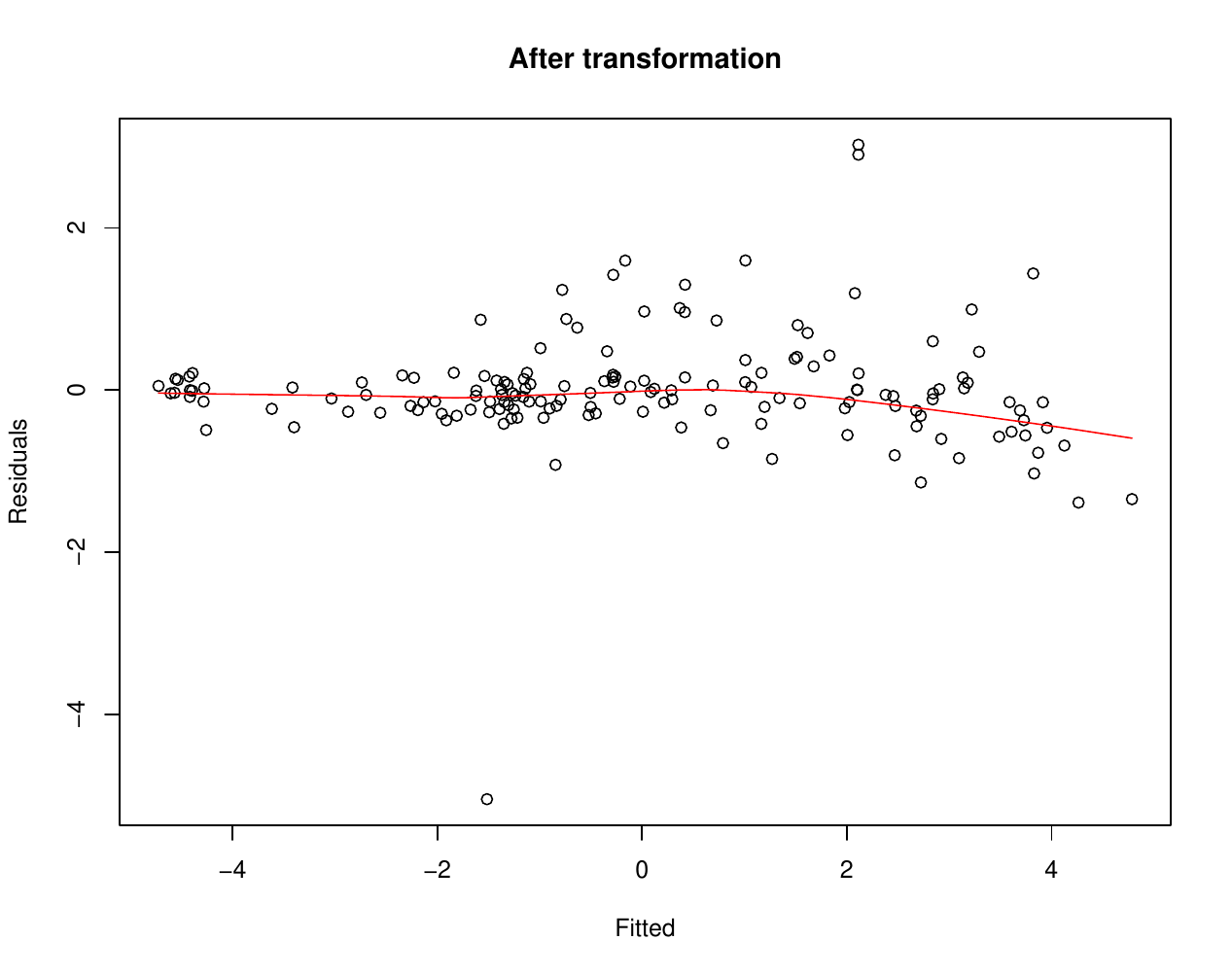}
        \caption{Plots of the residuals versus the centered heights and widths (top left and right) and of the residuals versus the fitted values (bottom). The added curve in red is a smoothing locally weighted polynomials obtained using the R funtion \text{lowess}.}
        \label{trans}
   \end{figure}

Next, we apply the transformation $y \mapsto y^{1/3}$ to the weights. We retain this time only the height and width which seem to have the most significant effects on the weight. Fitting a linear model after centering the transformed response and the 2  covariates yields $R^2 = 0.9133$, Also, the plots in Figure \ref{trans} show that the transformation substantially improves the joint behavior of the residuals and predictive variables.

\begin{table}[!h]
\begin{center}
\caption{The OLSE obtained for the fish market dataset after centering the transformed response and the retained continuous predictors height and width. The estimated standard deviations and  p-values associated with significance are also reported. }
\label{OLS}
\begin{tabular}{|ccc|} 
 \hline
  Predictor & Height   & Width  \\ [0.5ex] 
$\widehat \beta^{\text{OLSE}}_i$ & $0.099$  & $ 1.188$ \\
St. dev. & $0.022$  &   $0.056$      \\ 
p-value &  $1.37 \times 10^{-5} $ &  $<  2 \times 10^{-16}$   \\
 \hline
 \end{tabular}
\end{center}
\end{table}

\begin{figure}[H]
\centering
    \includegraphics[width=0.65\textwidth]{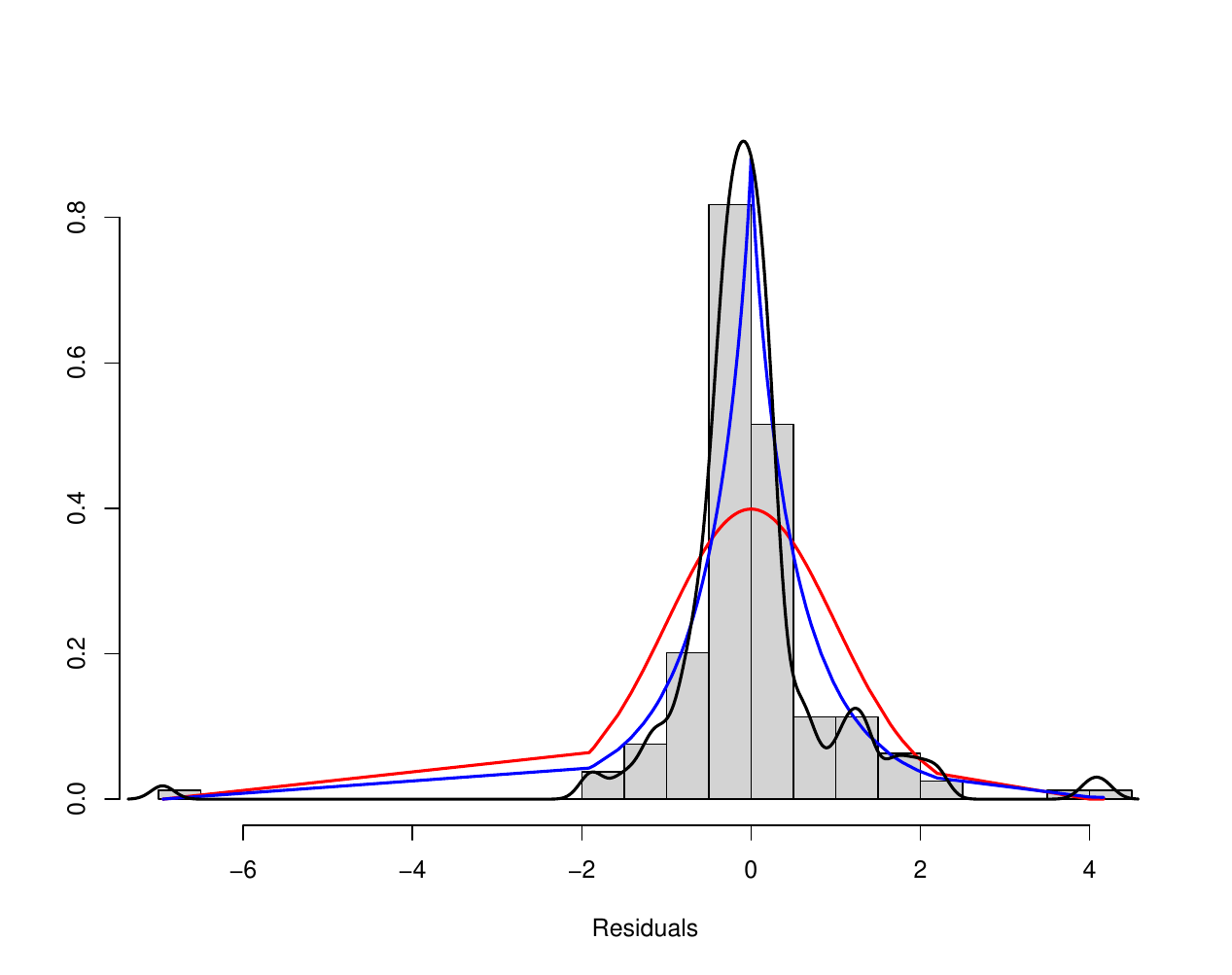}
        \caption{Histogram of the residuals of the regression fit. The curve in red/blue/black depicts the density of a standard Gaussian/the density $f_{0.85,1}$/a kernel density estimator.}
        \label{Histres}
   \end{figure}

Applying the Shapiro test for Gaussianity to the obtained residuals rejects the null hypothesis very significantly with p-value equal to $1.18 \times 10^{-14}$. The histogram of the rescaled residuals (with sample variance equal to 1) is shown in Figure \ref{Histres}. In the same figure, the density of a standard Gaussian in red clearly shows that the Gaussian assumption for the noise is indeed not appropriate.   The  curve in blue shows the density $f_{\gamma,1}$ with $\gamma = 0.85$. Finally, the curve in black depicts a standard kernel density estimator. Note that the rescaled residuals shows a slight shift to the left with a sample median equal to $-0.0824$. Although this suggests that one might find a more appropriate model for the noise distribution, we opted for simplicity especially that $f_{0.85,1}$ exhibits clearly a much better fit than the standard Gaussian density

The MLE of $\beta_0$ and $s_0$, $(\widehat \beta, \widehat s)$ is found by minimizing (-) the log-likelihood in $(\beta, s) \in \mathbb R^2 \times (0, \infty)$. The obtained estimators $\widehat \beta$ and $\widehat s$ and their estimated standard deviations  are shown in Table \ref{MLE}. The standard deviations are estimated using a resampling technique where we randomly picked 500 subsets of size 100. 
\begin{table}[!h]
\begin{center}
\caption{The MLE  of $(\beta_0, s_0)$ obtained for the fish market dataset after centering the transformed response (weight) and the continuous predictors height and width.}
\label{MLE}
\begin{tabular}{|cccc|} 
 \hline
 Predictor/scale & Height  & Width  & $s_0$ \\ [0.5ex] 
 MLE ($\widehat \beta^{\text{MLE}}_i / \widehat s^{\text{MLE}} $) & 0.132 &  1.106  & 0.569  
  \\  [0.5ex]
  St. dev. &  0.009 &    0.020 &   0.048  \\ [0.5ex]
 \hline
\end{tabular}
\end{center}
\end{table}

\medskip

To compare the predictive performances of the estimators, we take randomly $N$ data points as the training set while the remaining $159-N$ data points serve as a test set. This split is done for 500 replications in which the training set used to compute the OLSE and MLE and the test set to evaluate their respective mean square errors (MSE). Table \ref{mseare} gives the average value of the MSE of the OLSE and MLE as well as the estimated ARE $\widehat \eta$ for $N \in \{80, 90, 100 \}$. Note that to compute $\widehat \eta$ we used the expression
$$
\hat \eta =   \left(\frac{\det(\widehat{\text{Var}}
(\hat{\beta}_{\text{MLE}}))}{\det(\widehat{\text{Var}}
(\hat{\beta}_{\text{OLSE}}))}\right)^{1/2}
$$
where $\widehat{\text{Var}}(\hat{\beta}_{\text{MLE}})$  and $\widehat{\text{Var}}(\hat{\beta}_{\text{OLSE}})$  are the sample covariances of the MLE and OLSE obtained using the 500 replications. 

\begin{table}[!h]
\begin{center}
\caption{The average values of the MSE of the OLSE and MLE (in this order) and the sample ARE $\hat \eta$ computed using 500 replications when $N$, the size of the training subset, is $\in \{80, 90, 100 \}$.}
\label{mseare}     
\begin{tabular}{|ccc|} 
 \hline
 $N$ &  Av. MSE & $\widehat \eta$   \\ [0.5ex] 
 $80$  & (0.538, 0.531) & 0.266   \\ [0.5ex]
 $90$  &  (0.539, 0.532) &  0.280 \\  [0.5ex]
 $100$  & (0.528, 0.520) & 0.301   \\  [0.5ex]
\hline
\end{tabular}
\end{center}
\end{table}

Although the difference in the predictive performances is not very substantial the statistical gain of using the MLE over the OLSE for this data as reflected by the small values of $\widehat \eta  $ is not at all negligible. 

\section{Conclusions}
This paper addresses the problem of linear regression when the noise distribution is known up to a scale parameter. Unlike the Ordinary Least Squares Estimator (OLSE), which does not use such distributional knowledge, we considered the Maximum Likelihood Estimator (MLE) as a natural alternative and studied its theoretical and practical properties. 

We first proved that the likelihood function admits at least one maximizer under mild conditions, thereby guaranteeing existence of the MLE of both the regression vector $\beta_0$ and scale parameter and $s_0$. Second, we derived consistency of the MLE and then proved the joint asymptotic normality under suitable regularity conditions.  

 Through both theoretical derivations in Section \ref{sec: AsympEffic}, illustrated by simulation studies in Section \ref{simu}, we showed that the MLE is asymptotically at least as efficient as the OLSE and that the statistical gain can be substantial for non-Gaussian error distributions such as Weibull- or Gamma-based scale families. In these cases, we identified parameter regions where the MLE’s advantage is most pronounced, demonstrating that using knowledge of the noise distribution is far from being a marginal improvement, it can significantly reduce estimation variance in practice. 
 
The implications of our results are twofold. On the theoretical side, we provide a clear framework for understanding regression estimation when the error distribution is known up to a scale parameter, establishing existence, consistency, and asymptotic normality of the MLE. While these properties are classical in maximum likelihood theory, our contribution is to explicitly show how they apply in this specific linear regression setting and to quantify the relative efficiency compared to the OLSE. From a practical perspective, our results suggest that when the error distribution is reasonably well specified, the MLE can offer measurable improvements over the OLSE, particularly for non-Gaussian errors with heavier tails or higher peakedness. These improvements, while not surprising in light of likelihood theory, provide concrete guidance for applications in fields such as biomedical research or economics, where knowing the error structure up to a scale parameter is feasible. Building on the motivation of our work—following \citet{balabdaoui} in the context of semi-supervised learning where part of the data is unmatched and another part is matched; future research could explore practical applications of these results. In particular, it would be valuable to investigate how the efficiency gain of the MLE varies with the proportion of the dataset that is matched.

\section*{Appendix A: Notation and some definitions} 

\subsection*{Notation}
The notation and definitions presented in this section are taken from \cite{aadbookE2}.
\begin{itemize}
\item For a probability measure $\mathbb P$ defined on some measurable space $\mathcal X$, and integrable function $f$, $\mathbb P =  \int_{\mathcal X} f(x)  d\mathbb P(x) = \mathbb E_{\mathbb P}[f(X)]$ for a random variable $X \sim Q$. 

\item If $\mathbb P_n$ denotes the empirical probability measure associated with   i.i.d. random variables $X_1, \ldots, X_n \sim \mathbb P$, then we write
\begin{eqnarray*}
\mathbb G_n =  \sqrt n (\mathbb P_n -  \mathbb P),
\end{eqnarray*}
which is a signed random measure.  For an integrable function $f$, 
\begin{eqnarray*}
\mathbb G_n f & =   &  \sqrt (\mathbb P_n f - \mathbb P f)  =  \sqrt n \left( \frac{1}{n}  \sum_{i=1}^n f(X_i)  -  \mathbb E_{\mathbb P} E[f(X)]  \right)
\end{eqnarray*}
with $X \sim \mathbb P$.

\item  For some class of measure functions  $\mathcal F$, we write
$$
\Vert \mathbb G_n  \Vert_{\mathcal F}  = \sup_{f \in \mathcal F}  \vert \mathbb G_n f \vert.
$$

\end{itemize}

\subsection*{Covering numbers and entropy} 

\begin{definition}\label{Def1}

Let $\mathcal F$ be a class of measurable functions endowed by some norm $\Vert \cdot \Vert$. For a  given $\epsilon > 0$, the $\epsilon$-covering number $N(\epsilon, \mathcal F, \Vert \cdot \Vert)$ is defined as the minimal number of balls $\{f  \in \mathcal F: \Vert  f- f_i \vert < \epsilon \}$ needed to cover $\mathcal F$. Note that the centers $f_i$ do not need to belong to $\mathcal F$. The entropy is defined as $\log N(\epsilon, \mathcal F, \Vert \cdot \Vert)$. 
\end{definition}

\begin{definition}\label{Def2}

Let $\mathcal F$ be a class of measurable functions endowed by some norm $\Vert \cdot \Vert$. Suppose that $\mathcal F$ admits an envelope $F$; i.e., $\vert f \vert  \le  F$ for all $f \in \mathcal F$. Then, for a given $\delta > 0$ the uniform entropy integral is defined by
\begin{eqnarray*}
J(\delta, \mathcal F)  =  \sup_{Q} \int_0^\delta \sqrt{ 1 +  \log N(\epsilon \Vert F \Vert_{Q,2}, \mathcal F, L_2(Q) )} d\epsilon   
\end{eqnarray*}
where $\Vert F \Vert_{Q,2} =  \sqrt{\int F^2 dQ} $, and the supremum is taken over all probability measures with $\Vert F \Vert_{Q,2} > 0$.

\end{definition}

\section*{Appendix B: Proofs}

\subsection{Proofs for Section \ref{sec: MLE}}

\paragraph{Proof of Theorem \ref{exist}.} \ 
We will show that we can restrict attention to a compact set in $\mathbb R^d \times (0, \infty)$ when maximizing the log-likelihood. More precisely, we will exhibit an Euclidean ball $\mathcal{B}(0, R)$ and interval $[a, b]$ for some $R > 0, a > 0, b > 0$ (which depend on the data and $\alpha$) such that any $(\beta, s) \notin \mathcal{B}(0, R) \times [a, b]$ is not a good candidate for the maximization problem. First, note that 
\begin{eqnarray*}
 \frac{1}{n} \ell_n(\beta, s)  & =  &  -\log(s) + \frac{1}{n} \sum_{i=1}^n \log f \left( \frac{\vert Y_i-  \beta^T X_i\vert}{s} \right) \\
& \le & -\log(s)  + \log(\Vert f \Vert_\infty), 
\end{eqnarray*}
where $\Vert f\Vert_\infty  = \sup_{t \in \mathbb R} f(t) < \infty$  by Assumption (A3).  Hence, 
$$\lim_{s \to \infty} \sup_{\beta \in \mathbb R}  \ell_n(\beta, s) = -\infty.$$  
This implies that we must restrict the domain of $s$ to $(0, \sigma_1]$ for $\sigma_1 > 0$ which we will now exhibit. Since the maximum of $(\beta, s) \mapsto \ell_n(\beta, s)$ should be at least equal to $\ell_n(\beta_0, s_0)$ values of $s$ such that $\ell_n(\beta, s)  < \ell_n(\beta_0, s_0)$ for $\beta \in \mathbb R^d$ are not good candidates. Note that the inequality $\ell_n(\beta, s) <  \ell_n(\beta_0, s_0)$ holds if $s$ is such that
$$
-\log(s)  + \log(\Vert f \Vert_\infty) < \frac{1}{n} \ell_n(\beta_0, s_0) 
$$
or equivalently  $\displaystyle  s >  \exp\left( \log(\Vert f \Vert_\infty) - \frac{1}{n} \ell_n(\beta_0, s_0)  \right)$.  This means that good candidates for the maximization problem should be such that
\begin{eqnarray}\label{sigma1}
s \in (0,  \sigma_1], \ \text{where $\sigma_1: = \exp\left( \log(\Vert f \Vert_\infty) - \frac{1}{n} \ell_n(\beta_0, s_0)  \right)  $}. 
\end{eqnarray}
Next, we will show that good candidates for $\beta$ should have a bounded norm. Using the known inequality
\begin{eqnarray*}
\vert a + b \vert^\alpha \le c_\alpha (\vert a \vert^\alpha +  \vert b \vert^\alpha)
\end{eqnarray*}
with $c_\alpha = 1$ if $\alpha \in (0,1)$ and $2^{\alpha -1}$ if $\alpha \ge 1$, we can write that 
\begin{eqnarray*}
\vert Y_i - \beta^T X_i \vert^\alpha \ge c_\alpha^{-1} \vert \beta^T X_i \vert^\alpha -  \vert Y_i \vert^\alpha
\end{eqnarray*}
for $i =1, \ldots, n$ and hence 
\begin{eqnarray*}
\frac{1}{n} \ell_n(\beta, s) & \le  & \log(C)  - \log(s)  - 
 \frac{c_\alpha^{-1}}{s^\alpha} \frac{1}{n} \sum_{i=1}^n \vert \beta^T X_i \vert^\alpha +  \frac{1}{s^\alpha} \frac{1}{n} \sum_{i=1}^n  \vert Y_i \vert^\alpha  \\
 & = & \log(C)  - \log(s)  - 
 c_\alpha^{-1} \frac{\Vert \beta \Vert^\alpha}{s^\alpha} \frac{1}{n}\sum_{i=1}^n \vert {u_\beta}^T X_i \vert^\alpha +  \frac{1}{s^\alpha} \frac{1}{n} \sum_{i=1}^n  \vert Y_i \vert^\alpha  \\
\end{eqnarray*}
where $u_\beta \in \mathcal{S}_d$ the unit $(d-1)$-dimensional sphere.  Consider the function 
\begin{eqnarray*}
u \mapsto \frac{1}{n} \sum_{i=1}^n \vert u^T X_i \vert^\alpha    
\end{eqnarray*}
for $u \in \mathcal S_d$. Since this function is continuous on the compact set $\mathcal S_d$, it attains its minimum at some unit vector $u_*$. By assumption, we know that there exists at least an index $i \in \{1, \ldots, n \}$ such that $u_*^T X_i \ne 0$. This implies that $ \displaystyle \frac{1}{n} \sum_{i=1}^n \vert u^T_* X_i \vert^\alpha > 0$.   Thus,
 \begin{eqnarray*}
\frac{1}{n} \ell_n(\beta, s) & \le  & \log(C)  - \log(s)  - c_\alpha^{-1}    \frac{\Vert \beta \Vert^\alpha}{s^\alpha} \frac{1}{n} \sum_{i=1}^n \vert {u_*}^T X_i \vert^\alpha +  \frac{1}{s^\alpha} \frac{1}{n} \sum_{i=1}^n  \vert Y_i \vert^\alpha  \\
& = & \log(C)  - \log(s)  -  B_\alpha \frac{\Vert \beta \Vert^\alpha}{s^\alpha}  + \frac{A_\alpha}{s^\alpha}
\end{eqnarray*}
where 
\begin{eqnarray*}
A_\alpha =  \frac{1}{n} \sum_{i=1}^n  \vert Y_i \vert^\alpha, \ \ \text{and}  \  \ B_\alpha=  c_\alpha^{-1}  \frac{1}{n} \sum_{i=1}^n \vert {u_*}^T X_i \vert^\alpha.
\end{eqnarray*}
Writing $-\log(s) =   \log(1/s^\alpha)/\alpha$ and using the fact that $\log(t) \le t$ for all $t > 0$, it follows that
\begin{eqnarray*}
\frac{1}{n} \ell_n(\beta, s) & \le  & \log(C) + \left(\frac{1}{\alpha}  +  A_\alpha\right)  \frac{1}{s^\alpha}  -  B_\alpha \frac{\Vert \beta \Vert^\alpha}{s^\alpha}. 
\end{eqnarray*}
Then, $(\beta, s)$ is not a good candidate if $\ell_n(\beta, s) < \ell_n(\beta_0, s_0)$. Using the preceding inequality,  this holds if  
\begin{eqnarray}\label{ineq}
\log(C) + \left(\frac{1}{\alpha} + A_\alpha \right) \frac{1}{s^\alpha} - B_\alpha \frac{\Vert \beta \Vert^\alpha}{s^\alpha}  < \frac{1}{n} \ell_n(\beta_0, s_0).  
\end{eqnarray}
Suppose that $\displaystyle \log(C) \le \frac{1}{n} \ell_n(\beta_0, s_0)$, then the inequality in (\ref{ineq}) holds if 
$$
\left(\frac{1}{\alpha}  +  A_\alpha \right) \frac{1}{ s^\alpha} - B_\alpha \frac{\Vert \beta \Vert^\alpha}{s^\alpha} < 0
$$
or equivalently if
\begin{eqnarray*}
\Vert \beta \Vert > \left(\frac{1 + \alpha A_\alpha}{\alpha B_\alpha}\right)^{1/\alpha}.
\end{eqnarray*}
If $\displaystyle \log(C) > \frac{1}{n} \ell_n(\beta_0, s_0)$, then the inequality in (\ref{ineq}) holds if 
$$
\Vert \beta \Vert > \left \{ \frac{1}{B_\alpha} \left( \frac{1}{\alpha}  +  A_\alpha +  s^\alpha \left( \log(C) -   \frac{1}{n} \ell_n(\beta_0, s_0) \right)  \right) \right\}^{1/\alpha}
$$
which occurs if 
$$
\Vert \beta \Vert > \left \{ \frac{1}{B_\alpha} \left[ \frac{1 + \alpha A_\alpha}{\alpha} + \sigma_1^\alpha \left( \log(C) -   \frac{1}{n} \ell_n(\beta_0, s_0) \right)  \right] \right\}^{1/\alpha}
$$
where $\sigma_1$ is the same as in (\ref{sigma1}).   This means that good candidates from $\beta$ are those such that 
\begin{eqnarray*}
\Vert \beta \Vert \le 
\begin{cases}
\left(\frac{1 + \alpha A_\alpha}{\alpha B_\alpha}\right)^{1/\alpha}, \hspace{5.4cm} \text{if  $\displaystyle \log(C) \le \frac{1}{n} \ell_n(\beta_0, s_0)$}  \\
\left \{ \frac{1}{B_\alpha} \left[ \frac{1 + \alpha A_\alpha}{\alpha}  +  \sigma_1^\alpha \left( \log(C) -   \frac{1}{n} \ell_n(\beta_0, s_0) \right)  \right] \right\}^{1/\alpha}, \ \textrm{otherwise}.
\end{cases}
\end{eqnarray*}
We conclude that we should restrict attention to $\beta: \Vert \beta \Vert \le R_0$, where
\begin{eqnarray}\label{R0}
R_0 & = & \left(\frac{1 + \alpha A_\alpha}{\alpha B_\alpha}\right)^{1/\alpha} \vee \left \{ \frac{1}{B_\alpha} \left( \frac{1 + \alpha A_\alpha}{\alpha} +  \sigma_1^\alpha \left \vert \log(C) -   \frac{1}{n} \ell_n(\beta_0, s_0) \right\vert  \right) \right\}^{1/\alpha} \notag \\
& = &  \left \{ \frac{1}{B_\alpha} \left( \frac{1 + \alpha A_\alpha}{\alpha} +  \sigma_1^\alpha \left \vert \log(C) -   \frac{1}{n} \ell_n(\beta_0, s_0) \right\vert  \right) \right\}^{1/\alpha}.
\end{eqnarray}
Next, we will show that there exists $\sigma_0 > 0$ such that good candidates of $s$ should be at least equal to $\sigma_0$.  We have that 
\begin{eqnarray*}
\frac{1}{n}\ell_n(\beta, s) &\le  &  \log(C) - \log(s) - \frac{1}{s^\alpha} \sum_{i=1}^n \vert Y_i - \beta^T X_i \vert^\alpha.
\end{eqnarray*}
Consider the function 
$$
\beta \mapsto \sum_{i=1}^n \vert Y_i - \beta^T X_i \vert^\alpha
$$
for $\beta \in \mathcal{B}(0, R_0) = \{\beta \in \mathbb R^d: \Vert \beta  \Vert \le R_0 \}$. It follows from continuity of this function and compactness of $\mathcal{B}(0, R_0)$ that there exists $\beta_* \in \mathcal{B}(0, R_0)$ such that 
$$
\inf_{\beta \in \mathcal{B}(0, R_0)} \frac{1}{n} \sum_{i=1}^n \vert Y_i - \beta^T X_i \vert^\alpha  = \frac{1}{n} \sum_{i=1}^n \vert Y_i - \beta^T_* X_i \vert^\alpha.
$$
Furthermore, $\displaystyle  \frac{1}{n} \sum_{i=1}^n \vert Y_i - \beta^T_* X_i \vert^\alpha > 0$. In fact, if it were equal to $0$, then this would imply that for all $i=1, \ldots, n$  $Y_i - \beta^T_* X_i =0$, a case that is excluded by assumption.   Define now $\gamma = \alpha/2$.
We have that
\begin{eqnarray*}
\frac{1}{n}\ell_n(\beta, s) &\le  &  \log(C) - \log(s) - \frac{1}{s^\alpha} \frac{1}{n} \sum_{i=1}^n \vert Y_i - \beta^T_* X_i \vert^\alpha \\
& = &  \log(C) + \frac{1}{\gamma}\log(s^{-\gamma}) - \frac{1}{s^\alpha} \frac{1}{n} \sum_{i=1}^n \vert Y_i - \beta^T_* X_i \vert^\alpha  \\
& \le & \log(C) + \frac{2}{\alpha s^{\alpha/2}}  - \frac{1}{s^\alpha} \frac{1}{n} \sum_{i=1}^n \vert Y_i - \beta^T_* X_i \vert^\alpha.
\end{eqnarray*}
Then, the right side of the previous display is $ <  \frac{1}{n} \ell_n(\beta_0, s_0)$ if
\begin{eqnarray}\label{Ineq}
\log(C) + \frac{2}{\alpha s^{\alpha/2}} - \frac{1}{s^\alpha} \frac{1}{n} \sum_{i=1}^n \vert Y_i - \beta^T_* X_i \vert^\alpha < \frac{1}{n} \ell_n(\beta_0, s_0).
\end{eqnarray}
If $\displaystyle \log(C) \le \frac{1}{n} \ell_n(\beta_0, s_0)$, then the preceding inequality holds if
$$
s <  \left(\frac{\alpha}{2}\frac{1}{n} \sum_{i=1}^n \vert Y_i - \beta^T_* X_i \vert^\alpha\right)^{2/\alpha}.
$$
If $\displaystyle \log(C)  > \frac{1}{n} \ell_n(\beta_0, s_0)$, then put $t = s^{-\alpha/2}$ and $\tilde{A}_\alpha = n^{-1}  \sum_{i=1}^n \vert Y_i - \beta^T_* X_i \vert^\alpha $. The inequality in (\ref{Ineq}) is equivalent to 
\begin{eqnarray}\label{EqIneq}
t^2  \tilde{A}_\alpha - \frac{2}{\alpha} t + \frac{1}{n} \ell_n(\beta_0, s_0)  - \log(C)  > 0.
\end{eqnarray}
The discriminant of the second order polynomial on the left side is
$$
\Delta_\alpha =  \frac{4}{\alpha^2} - 4 \tilde{A}_\alpha \left(\frac{1}{n} \ell_n(\beta_0, s_0)  - \log(C) \right) > 0.
$$
The equation $t^2  \tilde{A}_\alpha - \frac{2}{\alpha} t + \frac{1}{n} \ell_n(\beta_0, s_0)  - \log(C) =0$ admits the roots
$$
t_1  =  \frac{1}{2 \tilde A_\alpha}\left( \frac{2}{\alpha} - \sqrt{\Delta_\alpha}\right), \ \ \text{and} \  \  t_2  =  \frac{1}{2 \tilde A_\alpha}\left( \frac{2}{\alpha} + \sqrt{\Delta_\alpha} \right).
$$
Note that $t_1 < 0 < t_2$. Since $t > 0$, the inequality in (\ref{EqIneq}) occurs if and only if 
$$
t > \frac{1}{2\tilde{A}_\alpha} \left( \frac{2}{\alpha} + \sqrt{\Delta_\alpha} \right)
$$
or equivalently 
$$
s < \left[\frac{1}{2\tilde{A}_\alpha} \left( \frac{2}{\alpha} + \sqrt{\Delta_\alpha} \right)\right]^{-2/\alpha}.
$$
We conclude that good candidates of $s$ should satisfy
\begin{eqnarray*}
s \ge
\begin{cases}
\left(\frac{\alpha}{2}\tilde{A}_\alpha\right)^{2/\alpha}, \ \ \hspace{2.7cm} \text{if $\displaystyle \log(C) \le \frac{1}{n} \ell_n(\beta_0, s_0)$ } \\
\left[\frac{1}{2\tilde{A}_\alpha} \left( \frac{2}{\alpha} + \sqrt{\Delta_\alpha} \right)\right]^{-2/\alpha}, \ \ \hspace{0.7cm} \text{otherwise}.
\end{cases}
\end{eqnarray*}
Define now 
\begin{eqnarray*}
\sigma_0 &=  & \left(\frac{\alpha}{2} \tilde A_\alpha \right)^{2/\alpha} \wedge \left[\frac{1}{2 \tilde A_\alpha} \left( \frac{2}{\alpha} + \sqrt{\widetilde{\Delta}_\alpha} \right)\right]^{-2/\alpha}, \\
&& \  \ \text{where} \ \widetilde{\Delta}_\alpha =   \frac{4}{\alpha^2} + 4 \tilde{A}_\alpha \left \vert \log(C) - \frac{1}{n} \ell_n(\beta_0, s_0)  \right \vert.
\end{eqnarray*}
But note that 
$$
\left[\frac{1}{2 \tilde A_\alpha} \left( \frac{2}{\alpha} + \sqrt{\widetilde{\Delta}_\alpha} \right)\right]^{-2/\alpha} >  \left(\frac{1}{\tilde{A}_\alpha \alpha}  \right)^{-2/\alpha} = \left(\alpha \tilde A_\alpha \right)^{2/\alpha}
$$
and therefore
\begin{eqnarray*}
\sigma_0 = \left(\frac{\alpha}{2} \tilde A_\alpha \right)^{2/\alpha}.
\end{eqnarray*}
It follows from the calculations above that with we can restrict attention to the sub-space  
\begin{eqnarray}\label{set}
\mathcal{B}(0, R_0) \times [a_0, a_1]
\end{eqnarray}
where $a_0 = \sigma_0 \wedge \sigma_1$ and $a_1 = \sigma_0 \vee \sigma_1$.   Continuity of $(\beta, s) \mapsto n^{-1} \ell_n(\beta)$ and compactness of the set in (\ref{set}) imply that the log-likelihood admits at least one maximizer. \hfill $\Box$

\paragraph{Proof of Proposition \ref{Determ}.}
First, we recall that
\begin{eqnarray*}
 \sigma_1 = \exp\left( \log(\Vert f \Vert_\infty) - \frac{1}{n} \ell_n(\beta_0, s_0)  \right).    
\end{eqnarray*}
By Assumption (A5), and by the Strong Law of large Numbers, 
$$
\frac{1}{n} \ell_n(\beta_0, s_0)  \to  -\log(s_0) + \mathbb E_{(\beta_0, s_0)} \log f \left(\frac{Y- \beta^T_0  X}{s_0} \right) < \infty
$$
almost surely.  Hence, with probability 1, and for $n$ large enough
$$
\frac{1}{n} \ell_n(\beta_0, s_0)  \ge  -\log(s_0) + \mathbb E_{(\beta_0, s_0)} \log f \left(\frac{Y- \beta^T_0  X}{s_0} \right) -1
$$
which implies that for $n$ large enough
\begin{eqnarray}\label{sigma*1}
 \sigma_1 \le \sigma^*_1: = \exp\left( \log(\Vert f \Vert_\infty s_0) - \mathbb E_{(\beta_0, s_0)} \log f \left(\frac{Y- \beta^T_0  X}{s_0} \right) + 1  \right).    
\end{eqnarray}
with probability 1. Recall $R_0$ from (\ref{R0}). Since it depends on $A_\alpha$ and $B_\alpha$ we will first find limits for the random quantities. We have that 
\begin{eqnarray*}
 A_\alpha & = & \frac{1}{n} \sum_{i=1}^n  \vert Y_i \vert^\alpha   \\
 & = & \frac{1}{n} \sum_{i=1}^n  \vert \epsilon_i + \beta_0^T X_i \vert^\alpha  \\
 & \le & c_\alpha \left( \frac{1}{n} \sum_{i=1}^n  \vert \epsilon_i  \vert^\alpha + \frac{1}{n} \sum_{i=1}^n  \vert \beta_0^T X_i \vert^\alpha \right) \\
 & \le &  c_\alpha \left( \frac{1}{n} \sum_{i=1}^n  \vert \epsilon_i  \vert^\alpha + \Vert \beta_0 \Vert^\alpha \frac{1}{n} \sum_{i=1}^n  \Vert X_i \Vert^\alpha \right)
 \end{eqnarray*}
which converges to $c_\alpha \left( \mathbb E[\vert \epsilon\vert^\alpha] +  \Vert \beta_0 \Vert^\alpha \mathbb E[\Vert X\Vert^\alpha]  \right)  < \infty
$ almost surely by the SLLN, which can be applied by Assumption (A6) and the fact that $\displaystyle \mathbb E[\vert \epsilon \vert^\alpha ] \le Cs^{-1}_0 \int_{\mathbb R} \vert t\vert^\alpha \exp\left(- \vert t /s_0 \vert^\alpha\right) dt < \infty$.  This implies that almost surely
\begin{eqnarray}\label{A*}
 A_\alpha \le 2 c_\alpha \left( \mathbb E[\vert \epsilon\vert^\alpha] +  \Vert \beta_0 \Vert^\alpha \mathbb E[\Vert X\Vert^\alpha]  \right)  := A^*
\end{eqnarray}
for $n$ large enough.

Also,
\begin{eqnarray*}
B_\alpha & =   & c^{-1}_\alpha  \frac{1}{n} \sum_{i=1}^n \vert u^T_* X_i  \vert^\alpha \\
& =  & c^{-1}_\alpha  \inf_{u \in \mathcal S_{d}} \frac{1}{n} \sum_{i=1}^n \vert u^T_* X_i \vert^\alpha \\
& \rightarrow & c^{-1}_\alpha  \inf_{u \in \mathcal S_{d}} \mathbb E[\vert u^T X \vert^\alpha]
\end{eqnarray*}
almost surely using the SLLN and the continuity Theorem. Hence, with probability 1
$$
B_\alpha \ge \frac{1}{2} c^{-1}_\alpha \inf_{u \in \mathcal S_{d}} \mathbb E[\vert u^T X \vert^\alpha]: = B^*
$$
for $n$ large enough. Note that $\inf_{u \in \mathcal S_{d}} \mathbb E[\vert u^T X \vert^\alpha] > 0$. In fact, by compactness of $\mathcal S_d$ and continuity 
of the function $u \mapsto \mathbb E[\vert u^T X \vert^\alpha]$ (this can be easily shown using the Dominance Convergence theorem), there exists $v \in \mathcal{S}_d$ such that $\inf_{u \in \mathcal S_{d}} \mathbb E[\vert u^T X \vert^\alpha] = \mathbb E[\vert v^T X \vert^\alpha]$. If the latter is equal to $0$, then this would imply that $\mathbb P(v^T X = 0) =1$ and hence $v =0$ by Assumption (A2) which implies that $v^T X$ must have an absolutely continuous distribution in case $v \ne 0$. Finally, 
\begin{eqnarray*}
\left \vert  \log(C)  - \frac{1}{n} \ell_n(\beta_0, s_0) \right \vert  & \le  &    \vert \log(C)  \vert +   \left \vert \frac{1}{n} \ell_n(\beta_0, s_0)  \right \vert  \\
& \le &  \vert \log(C)  \vert  + \vert \log(s_0) \vert + \left \vert \mathbb E_{(\beta_0, s_0)} \log f \left(\frac{Y- \beta^T_0  X}{s_0} \right)  \right \vert +1 \\
&: = & \Delta^*
\end{eqnarray*}
with probability 1 for $n$ large enough. We conclude that 
\begin{eqnarray}\label{R*}
R_0 \le \left \{ \frac{1}{B^*}  \left( 
\frac{1}{\alpha} + A^* +  \sigma^*_1  \Delta^* \right) \right \}^{1/\alpha}: = R^* 
\end{eqnarray}
with probability 1 for $n$ large enough.  Finally,
recall that
\begin{eqnarray*}
\sigma_0 =  \left(\frac{\alpha \tilde{A}_\alpha}{2}  \right)^{2/\alpha}  
\end{eqnarray*}
where
\begin{eqnarray*}
 \tilde{A}_\alpha &= &  \frac{1}{n} \sum_{i=1}^n \sum_{i=1}^n \vert Y_i - \beta^* X_i \vert^\alpha  \\
 & \le & c_\alpha \left( \frac{1}{n} \sum_{i=1}^n \sum_{i=1}^n \vert Y_i \vert^\alpha +  \Vert \beta^* \Vert^\alpha \frac{1}{n} \sum_{i=1}^n \Vert X_i \Vert^\alpha  \right) \\
 & \le  &  2 c_\alpha \left( \mathbb E[\vert \epsilon\vert^\alpha] +  \Vert \beta_0 \Vert^\alpha \mathbb E[\Vert X\Vert^\alpha]  \right)  + 2 c_\alpha  R^* \mathbb E[\Vert X \Vert^\alpha] \\
 & = & A^* + 2 c_\alpha  R^* \mathbb E[\Vert X \Vert^\alpha] := \tilde{A}^*_\alpha, \ \ \text{where $A^*$ is the same as in (\ref{A*}) }
\end{eqnarray*}
with probability 1 for $n$ large enough. Thus,
\begin{eqnarray}\label{sigma*0}
\sigma_0 \le  \left(\frac{\alpha \tilde{A}^*_\alpha}{2}\right)^{2/\alpha}:= \sigma^*_0.
\end{eqnarray}
Using the calculations above, we conclude that with probability 1 there exists $n_0$ such that for all $n \ge n_0$ the likelihood maximization problem can be restricted on the compact set 
$$
\mathcal{B}(0, R^*) \times [a^*, b^*]
$$
where $R^*$ is the radius defined in (\ref{R*}),  $a^* = \sigma^*_0 \wedge \sigma^*_1$ and $b^* = \sigma^*_0 \vee \sigma^*_1$, with $\sigma^*_0$ and $\sigma^*_1$ are the same as in (\ref{sigma*0}) and  (\ref{sigma*1}) respectively.  \hfill $\Box$

\paragraph{Proof of Theorem \ref{Consis}.}
We will start with showing the following properties:
\begin{eqnarray}\label{UC}
\sup_{(\beta, s) \in \mathcal{B}(0, R^*) \times [a^*, b^*]}  \left \vert \frac{1}{n} \ell_n(\beta, s) - \ell(\beta, s) \right \vert = o_{\mathbb P}(1)
\end{eqnarray}
with
\begin{eqnarray*}
\ell(\beta, s) &= & \mathbb E_{(\beta_0, s_0)} \log \left[\frac{1}{s} f \left( \frac{Y - \beta^T X}{s}\right) \right]  \\
&= & \int \log \left[\frac{1}{s} f \left( \frac{y - \beta^T y}{s}\right) \right]  \frac{1}{s_0} f \left( \frac{y - \beta^T_0 x}{s_0}\right) f^X(x) dx dy,
\end{eqnarray*}
and for any $r_1, r_2 > 0$ 
\begin{eqnarray}\label{Iden}
\sup_{\beta \in \mathcal O_{r_1, r_2}} \ell(\beta, s) <   \ell(\beta_0, s_0)
\end{eqnarray}
where 
$$
\mathcal O_{r_1, r_2} = \{(\beta, s):  \Vert \beta - \beta_0 \Vert > r_1, \vert s - s_0 \vert > r_2 \}.
$$
The first result is known as uniform consistency while the second is the identifiability condition.   To show (\ref{UC}) we will resort to some tools from empirical process theory.  we refer to Appendix A for the notation that is commonly used in this scope.  We can write that
\begin{eqnarray*}
\frac{1}{n} \ell_n(\beta, s)  -  \ell(\beta, s)  & =  &  \int \log \left( \frac{1}{s}  f\left( \frac{y-\beta^T x}{s}  \right)   \right) d(\mathbb P_n -  \mathbb P)(x, y)  \\
& = &  \int \log \left( f\left( \frac{y-\beta^T x}{s}  \right)   \right) d(\mathbb P_n -  \mathbb P)(x, y)
\end{eqnarray*}
where $\mathbb P$ and $\mathbb P_n$ denote the true probability measures associated with the distribution of $(X, Y)$ and empirical probability measure based on the random sample $(X_1, Y_1), \ldots, (X_n, Y_n)$. Note that the second equality follows from the fact that $ \int d(\mathbb P_n - \mathbb P)(x, y) =0$. 

Consider now the class of functions
\begin{eqnarray*}
 && \{(x, y) \mapsto  l_{\beta, s}(x, y) = s^{-1}(y - \beta^T x), \ \beta \in \mathcal{B}(0, R^*), \ s \in [a^*, b^*] \}  \\ 
 & & \subset   \{(x, y) \mapsto  l_{\beta}(x, y) = y - \beta^T x, \ \beta \in \mathcal{B}(0, R^*/a^*), \ s \in [a^*, b^*] \} \\
 && := \mathcal{L}. 
\end{eqnarray*}
In fact $s^{-1}( y - \beta^T x)  =  y' - \beta'^T x$ where $y' = s^{-1} y \in \mathbb R$ and $\vert \beta' \vert =  s^{-1} \Vert \beta \Vert \le R^*/a^*$. 

The class $\mathcal{L}$ is indexed by $\beta$ and hence is a finite dimensional vector space. From \cite[Lemma 2.6.16]{aadbookE2} it follows that $\mathcal L$ is a VC-subgraph of dimension $V \le d+2$. On the other hand, we can write that 
\begin{eqnarray*}
\log f(t) & = &  \log f(t) \mathds{1}_{t \ge 0}  +   \log f(t) \mathds{1}_{t < 0}  \\
& =  &  m_+(t)  + m_{-}(t)  -  \log(f(0))
\end{eqnarray*}
where $m_+(t)  =  \log f(t) \mathds{1}_{t \ge 0} + \log(f(0)) \mathds{1}_{t < 0} $ and $m_{-}(t) = \log f(t) \mathds{1}_{t < 0} + \log(f(0)) \mathds{1}_{t \ge 0}$.  The functions $m_+$ and $m_-$ are monotone non-increasing and non-decreasing respectively.  Then, using again the fact that $\int d(\mathbb P_n - \mathbb P) =0$ we can write for any $\beta \in \mathcal{B}(0, R^*/a^*) \times [a^*, b^*]$
\begin{eqnarray*}
\int \log f \circ l_{\beta}(x, y) d(\mathbb P_n - \mathbb P)(x, y)  & = &   \int m_-  \circ l_{\beta, s}(x, y)  d(\mathbb P_n - \mathbb P)(x, y)   \\
&& \  \  +   \int m_+  \circ l_{\beta}(x, y)  d(\mathbb P_n - \mathbb P)(x, y) 
\end{eqnarray*}
and hence
\begin{eqnarray*}
&& \sup_{(\beta, s) \in \mathcal{B}(0, R^*) \times [a^*, b^*]}  \left \vert \frac{1}{n} \ell_n(\beta, s) - \ell(\beta, s) \right \vert  \le  \\
&&  =   \sup_{\beta \in \mathcal{B}(0, R^*/a ^*)}  \left \vert \int m_-  \circ l_{\beta}(x, y)  d(\mathbb P_n - \mathbb P)(x, y)   \right \vert   \\
&&  \  \  +  \sup_{\beta \in \mathcal{B}(0, R^*/a^*)}  \left \vert  \int m_+  \circ l_{\beta}(x, y)  d(\mathbb P_n - \mathbb P)(x, y)  \right \vert.   
\end{eqnarray*}
By  \cite[Lemma 2.6.20 - (viii)]{aadbookE2} the class of functions $\mathcal G_-  =  m_{-}  \circ \mathcal L $
is a VC-subgraph with some finite index  $V_- > 0$.   Furthermore,  for $(\beta, s) \in \mathcal{B}(0, R^*) \times [a^*, b^*]$ and $(x, y) \in \mathbb R^d \times \mathbb R$  we have by the Cauchy-Schwarz inequality that $ s^{-1} \vert y - \beta x  \Vert \le a^{-1}_* (\vert y \vert +  R^* \Vert x \Vert) $ and the fact that $f$ is non-increasing on $[0, \infty)$
\begin{eqnarray*}
\log f\left(\frac{\vert y \vert + R^* \Vert x \Vert }{a^*}  \right)  \le  \log  f\left(\frac{y - \beta^T x}{s}\right) = \log f\left( \frac{\vert y - \beta^T x \vert}{s}\right )   \le \log f(0)
\end{eqnarray*}
and hence
\begin{eqnarray}\label{env}
\left \vert  \log  \left(\frac{y - \beta^T x}{s}\right)  \right \vert  & \le  &   \vert \log f(0) \vert \vee \left \vert \log f\left(\frac{\vert y \vert + R^* \Vert x \Vert }{a^*}  \right)    \right \vert \notag \\
& \le & \vert \log f(0) \vert +  \left \vert \log f\left(\frac{\vert y \vert + R^* \Vert x \Vert }{a^*}  \right)    \right \vert : = G(x, y).
\end{eqnarray}
By Theorem \cite[Lemma 2.6.7]{aadbookE2} we known that
\begin{eqnarray}
 N(\delta  \Vert G \Vert_{Q,2}, \mathcal{G}_-, L_2(Q) ) \le K \left(\frac{1}{\delta}\right)^{2 V_-}  
\end{eqnarray}
where $K > 0$ depends only on $V_-$, and which can be taken without loss of generality to be $> 1/e$. Above,  $Q$ is any probability measure such that $ \Vert G \Vert_{Q,2} = \sqrt{ \int G^2 dQ} > 0$. 

For $\eta > 0$, recall the the uniform entropy integral $J(\eta, \mathcal G_-)$ defined in Definition \ref{Def2}.
We have that 
\begin{eqnarray*}
J(1, \mathcal{G}_+) & \le  & \int_0^1  \left(\sqrt{1 + \log K} +  \sqrt{2 V_-} \sqrt{\log\left(\frac{1}{\delta} \right)} \right) d\delta  \\
& \le &  \sqrt{1 + \log K}  +  \sqrt{2 V_-} \int_0^1 \frac{1}{\sqrt{\delta}} d \delta \\
& =  & \sqrt{1 + \log K}  +  2\sqrt{2 V_-}.
\end{eqnarray*}
Also, Assumption (A8) implies that $\mathbb E_{(\beta_0, s_0)} G^2(X,Y)  = \Vert G \Vert^2_{L_2(\mathbb P)} < \infty.$  Hence,  by \cite[Theorem 2.14.1]{aadbookE2}, it follows that
\begin{eqnarray*}
\mathbb E[\Vert \mathbb G_n \Vert^2_{\mathcal{G}_-}]^{1/2} \lesssim J(1, \mathcal{G}_-) \times  \Vert  G \Vert_{L_2(\mathbb P)} < \infty. 
\end{eqnarray*}
Using Markov's inequality, we conclude that $\Vert \mathbb G_n \Vert_{\mathcal{G}_-}  = O_{\mathbb P}(1) $, and hence
\begin{eqnarray*}
\sup_{\beta \in \mathcal{B}(0, R^*/a ^*)}  \left \vert \int m_-  \circ l_{\beta}(x, y)  d(\mathbb P_n - \mathbb P)(x, y)   \right \vert  = O_{\mathbb P}(1/\sqrt n).
\end{eqnarray*}
Since a similar reasoning can be applied to the class $\mathcal{G}_+  = m_+ \circ \mathcal L$, it follows that 
$$
\sup_{(\beta, s) \in \mathcal{B}(0, R^*) \times [a^*, b^*]}  \left \vert \frac{1}{n} \ell_n(\beta, s) - \ell(\beta, s) \right \vert=  O_{\mathbb P}(1/\sqrt n)
$$
which in turn implies the uniform consistency in (\ref{UC}).  \\

\bigskip
Now, we show the identifiability condition in (\ref{Iden}). We first show that if $(\beta, s) \ne (\beta_0, s_0)$ then $\ell(\beta, s) < \ell(\beta_0, s_0)$. For any $(\beta, s) \in \mathbb R^d \times (0, \infty$ we have by Jensen's inequality applied to the strictly convex function $-\log $ that
\begin{eqnarray*}
\mathbb E_{(\beta_0, s_0)} \log \left[ \frac{1}{s} f\left( \frac{Y- \beta^T X}{s} \right) \right]  \le  \mathbb E_{(\beta_0, s_0)} \log \left[ \frac{1}{s_0} f\left( \frac{Y- \beta^T_0 X}{s_0} \right) \right]  
\end{eqnarray*}
with equality if and only if there exits a constant $C > 0$ such that
\begin{eqnarray*}
\mathbb P\left(\frac{1}{s} f\left( \frac{Y- \beta^T X}{s} \right) =  C  \frac{1}{s_0} f\left( \frac{Y- \beta^T_0 X}{s_0} \right) \right) =1
\end{eqnarray*}
Since both sides of the equality are densities, we must have $C=1$. Thus, $\ell(\beta, s) = \ell(\beta_0, s_0)$ if and only if 
\begin{eqnarray*}
\frac{1}{s} f\left( \frac{y- \beta^T x}{s} \right)  = \frac{1}{s_0} f\left( \frac{y- \beta^T_0 x}{s_0} \right)
\end{eqnarray*}
for $\mathbb P$-almost every $(x, y)$.  Fix $x$ and put $t = (y- \beta_0^T x)/s_0$. Also, let $c = s^{-1} (\beta_0 - \beta)^T x$ and $\lambda = s_0/s$. Then, the previous equality yields
\begin{eqnarray*}
\lambda  f( \lambda t  +c )  = f(t)
\end{eqnarray*}
for almost all $t$.   We show that $c=0$ and $\lambda =1$.  Using the identity above recursively, it follows that
\begin{eqnarray*}
f(t)  & = & \lambda^2  f( \lambda^2 t  + c (\lambda +1)) \\
& = &  \lambda^3 f( \lambda^3 t  + c (\lambda^2 +\lambda +1)) \\
& \vdots  &  \\
& = &  \lambda^k f( \lambda^k t  + c (\lambda^k +\lambda^{k-1} + \ldots +1))
\end{eqnarray*}
for all $k \ge 1$.  Hence, for $t = 0$, we have that
\begin{eqnarray}\label{f0}
f(0)  =  \lambda^k f \big( c (\lambda^k +\lambda^{k-1} + \ldots +1) \big)    
\end{eqnarray}
for all $k \ge 1$. If $c \ne 0$, then  this implies that
\begin{eqnarray*}
f(0)  =  \frac{\lambda^k}{c \sum_{j=0}^k \lambda^j }  c \sum_{j=0}^k \lambda^j \  f \big( c \sum_{j=0}^k \lambda^j \big)
\end{eqnarray*}
Since $f$ admits a finite expectation,  even and assumed to be continuous by Assumption (A3), we must have that $\lim_{\vert x \vert \to \infty}  \vert x \vert f(\vert x  \vert )  =0$.  If $\lambda \ne 1$, then by letting $k \to \infty$, we obtain that $f(0) =0$, which is impossible. If $\lambda = 1$, then the equality in (\ref{f0}) becomes $f(0)  = f( c( k+1)) $. By letting $k \to \infty$ we obtain again that $f(0) =0$. Hence, we must have $c = 0$. Then, this yields the identity $f(t) = \lambda^k f(\lambda^k t)$ and hence $f(0)  = \lambda^k f(0)$ for all $k \ge 1$, which implies that $\lambda =1$.  We conclude that $s = s_0$ and $(\beta - \beta_0)^T x =0$ for almost all $x$. By Assumption (A2), the distribution of $X$ is absolutely continuous and hence we must have $\beta - \beta_0 = 0$.  

Now, let $r_1 > 0, r_2  > 0$ and the set $\mathcal O_{r_1, r_2} = \{(\beta, s): \Vert \beta - \beta \vert > r_2, \vert s - s_0 \vert > r_2 \}$. Since we know that we can restrict ourselves to the compact set $C^* : = \{ \beta: \Vert \beta \Vert \le R^*  \times [a^*, b^*] \}$, then showing (\ref{Iden}) is equivalent to show that
$$
\sup_{ \beta \in \mathcal O_{r_1, r_2} \times C^*}  \ell(\beta, s) < \ell(\beta_0, s_0).
$$
Note first that 
$$
\sup_{ \beta \in \mathcal O_{r_1, r_2} \times C^*}  \ell(\beta, s) \le \sup_{ \beta \in \overline{\mathcal O}_{r_1, r_2} \cap C^*}  \ell(\beta, s)
$$
with $\overline{\mathcal O}_{r_1, r_2}$ the closure of $\mathcal O_{r_1, r_2}$.  Suppose that $\sup_{ \beta \in \overline{\mathcal O}_{r_1, r_2} \cap C^*}  \ell(\beta, s)  =  \ell(\beta_0, s_0)$.  By compactness of $ \overline{\mathcal O}_{r_1, r_2} \cap C^* $ and continuity of $(\beta, s) \mapsto \ell(\beta, s)$, this implies that there exits $(\beta^*, s^*)$ such that
 $$
\ell(\beta^*, s^*) =  \sup_{ \beta \in \overline{\mathcal O}_{r_1, r_2} \cap C^*}  = \ell(\beta_0, s_0).
$$
By the proof above, this means that $\beta^* = \beta_0$ and $s^* = s_0$, which is impossible. This completes the proof since by \cite[Corollary 3.2.3]{aadbookE2}, it follows that 
$$
(\widehat \beta_n, \widehat s_n)  \to_{\mathbb P} (\beta_0, s_0).
$$
\hfill $\Box$ 

\paragraph{Proof of Theorem \ref{asympMLE}.}  Using (A9), (A10) and Taylor expansion of $(\beta, s) \mapsto (1/n) \ell_n(\beta, s)$ we can write  that
\begin{eqnarray}\label{EqCTL}
\mathds{O}_{d+1}= 
 \left( 
\begin{array}{c}
\frac{1}{n} \sum_{i=1}^n \frac{\partial \log f_s(Y_i - \beta^T X_i)}{\partial \beta}|_{\theta = \theta_0} \\
\frac{1}{n} \sum_{i=1}^n \frac{\partial \log f_s(Y_i - \beta^T X_i)}{\partial s}|_{\theta = \theta_0}  
\end{array}
\right)  +  M_n \left(
\begin{array}{c}
\widehat \beta_n - \beta_0 \\
\widehat s_n - s_0
\end{array}
\right)
\end{eqnarray}
where $M_n$ is the $(d+1) \times (d+1)$ matrix given by 
\begin{eqnarray*}
M_n =  \left(
\begin{array}{cc}
 \frac{1}{n} \sum_{j=1}^n \frac{\partial^2 \log f_s(Y_i - \beta^T X_i)}{\partial \beta \partial \beta^T}|_{\theta = (\widetilde{\beta}_n, \tilde{s}_n)}   &   \frac{1}{n} \sum_{j=1}^n \frac{\partial^2 \log f_s(Y_i - \beta^T X_i)}{\partial \beta \partial s}|_{\theta = (\widetilde{\beta}_n, \tilde{s}_n)}  \\
   &   \\
 \frac{1}{n} \sum_{j=1}^n \frac{\partial^2 \log f_s(Y_i - \beta^T X_i)}{\partial \beta \partial s}|_{\theta = (\widetilde{\beta}_n, \tilde{s}_n)}  &  \frac{1}{n} \sum_{j=1}^n \frac{\partial^2 \log f_s(Y_i - \beta^T X_i)}{\partial s^2}|_{\theta = (\widetilde{\beta}_n, \tilde{s}_n)}
\end{array}
\right)
\end{eqnarray*}
with $(\widetilde{\beta}_n, \tilde{s}_n)$ belongs to a small neighborhood of $\theta_0$, a consequence of consistency of the MLE $\widehat \theta_n$; see Theorem \ref{Consis}.   From Assumptions (A9-A11) and the SLLN, it follows that
\begin{eqnarray*}
M_n =   M_0  + o_{\mathbb P}(1)  
\end{eqnarray*}
with 
\begin{eqnarray*}
M_0 =  \left(  
\begin{array}{cc}
\mathbb E_{\theta_0} \left[ \frac{\partial^2 \log f_s(Y - \beta^T X)}{\partial \beta \partial \beta^T}|_{\theta = \theta_0} \right]  &   \mathbb E_{\theta_0} \left[ \frac{\partial^2 \log f_s(Y - \beta^T X)}{\partial \beta \partial s}|_{\theta = \theta_0} \right] \\
  &  \\
\mathbb E_{\theta_0} \left[ \frac{\partial^2 \log f_s(Y - \beta^T X)}{\partial \beta \partial s}|_{\theta = \theta_0} \right]   &  \mathbb E_{\theta_0} \left[ \frac{\partial^2 \log f_s(Y - \beta^T X)}{\partial s^2}|_{\theta = \theta_0} \right]
\end{array}
\right).
\end{eqnarray*}
Next, we compute explicitly the entries of $M_0$.  We compute 
\begin{eqnarray*}
\frac{\partial^2 \log f_s(y - \beta^T x)}{\partial \beta \partial \beta^T }  & =  &   \frac{\partial^2}{\partial \beta \partial \beta^T} \left( \log(1/s)  +  \log f\left(\frac{y-\beta^T x }{s}\right)  \right)  \\
& = &  \frac{\partial }{\partial \beta^T} \left(  - \frac{1}{s}  \frac{f'\left(\frac{y-\beta^T x }{s}\right)}{f\left(\frac{y-\beta^T x }{s}\right)}  x\right)  \\
& = &    \frac{1}{s^2} \frac{f''(\frac{y-\beta^T x}{s})f(\frac{y-\beta^T x}{s})-\left(f'(\frac{y-\beta^T x}{s})\right)^2}{f(\frac{y-\beta^Tx}{s})^2}  \cdot xx^T. 
\end{eqnarray*}
Then, 
\begin{eqnarray*}
&& -\mathbb E_{\theta_0} \left[ \frac{\partial^2 \log f_s(Y - \beta^T X)}{\partial \beta \partial \beta^T}|_{\theta = \theta_0} \right]   \\
&& = - \frac{1}{s_0^2} \int  \frac{f''\left(\frac{y-\beta^T_0 x}{s_0}\right)f\left(\frac{y-\beta^T_0 x}{s}\right)-\left(f'\left(\frac{y-\beta^T_0 x}{s_0}\right)\right)^2}{f\left(\frac{y-\beta^T_0x}{s_0}\right)^2}  \cdot xx^T  \frac{1}{s_0} f\left(\frac{y-\beta_0^T x}{s_0}\right) f^X(x)  dx dy  \\
&&  = - \frac{1}{s_0^2} \int \left( f''(t) - \frac{(f'(t))^2}{f(t)} \right) dt \cdot \mathbb E(XX^T) \\
&& =  \frac{1}{s_0^2} \left(\int \frac{(f'(t))^2}{f(t)} dt  \right)   \cdot \mathbb E(XX^T)  = c_1 \mathbb E(XX^T).
\end{eqnarray*}
using the change of variable $t =  (y - \beta_0^T x)/s_0$ (for a fixed $x$), the fact that $\int f''(t) dt = 0$ and $c_1 \in \mathbb R$. To see the latter, note that for any $A > 0$
\begin{eqnarray*}
 \int_A^\infty \frac{(f'(t))^2}{f(t)}  dt  & = &  \int_A^\infty \frac{1}{t^2}  t^2 \frac{(f'(t))^2}{f(t)}  dt  \\
 & \le &  \frac{1}{A^2}  \int_A^\infty t^2 \frac{(f'(t))^2}{f(t)}  dt   \to 0
\end{eqnarray*}
as $A \to \infty$ using (A11).  Also,
\begin{eqnarray*}
&&\frac{\partial^2 \log f_s(y - \beta^T x)}{\partial \beta \partial s}  \\
&& =     \frac{\partial}{\partial \beta} \left( -\frac{1}{s} - \frac{f'\left(\frac{y-\beta^T x}{s}\right)}{f\left(\frac{y-\beta^T x}{s}\right)} \cdot \left(\frac{y-\beta^T x}{s^2}\right)  \right) \\
&& =  -\frac{f''\left(\frac{y-\beta^T x}{s}\right)f\left(\frac{y-\beta^T x}{s}\right)\left(-\frac{x}{s}\right)-\left(f'\left(\frac{y-\beta^Tx}{s}\right)\right)^2\left(-\frac{1}{s} x\right) }{f\left(\frac{y-\beta^T x}{s}\right)^2 } \cdot \frac{y-\beta^T x}{s^2} \\
&& \ - \ \frac{x}{s^2} \frac{f'\left(\frac{y-\beta^T x}{s}\right)}{f\left(\frac{y-\beta^T x}{s}\right)}. 
\end{eqnarray*}
Hence,
\begin{eqnarray*}
&& -\mathbb E_{\theta_0} \left[ \frac{\partial^2 \log f_s(Y - \beta^T X)}{\partial \beta \partial s}|_{\theta = \theta_0} \right]   \\
&& = \int  x\left(\frac{-f''\left(\frac{y-\beta^T_0 x}{s_0}\right)f\left(\frac{y-\beta^T_0 x}{s}\right) +\left(f'\left(\frac{y-\beta^T_0x}{s_0}\right)\right)^2}{f\left(\frac{y-\beta^T_0 x}{s_0}\right)^2 } \cdot \frac{y-\beta^T_0 x}{s^3_0}  \right) \frac{1}{s_0} f\left(\frac{y-\beta_0^T x}{s_0} \right) f^X(x) dx dy \\
&&  \ + \ \int x \frac{1}{s^2_0} \frac{f'\left(\frac{y-\beta^T x}{s}\right)}{f\left(\frac{y-\beta^T x}{s}\right)}  \frac{1}{s_0} f\left(\frac{y-\beta_0^T x}{s_0} \right) f^X(x) dx dy  \\
&& =  \frac{1}{s_0^2} \int \frac{t (f'(t))^2}{f(t)} dt \ \cdot \mathbb E(X)   + \frac{1}{s_0^2} \int f'(t) dt  \ \cdot \mathbb E(X)  = 0
\end{eqnarray*}
using the same change of variable as above, $\mathbb E[X] < \infty $ by Assumption (A6), that $\int f'(t) dt = \int f''(t) dt = 0$,  the fact that $ \int \vert t \vert  (f'(t))^2 /f(t) dt < \infty$ (implied by (A11)) and using that the function $t \mapsto t (f'(t))^2 /f(t)$ is odd.  Finally, we have that 
 \begin{eqnarray*}
&&\frac{\partial^2 \log f_s(y - \beta^T x)}{\partial s^2}  \\
&& =     \frac{\partial}{\partial s} \left( -\frac{1}{s} - \frac{f'\left(\frac{y-\beta^T x}{s}\right)}{f\left(\frac{y-\beta^T x}{s}\right)} \cdot \left(\frac{y-\beta^T x}{s^2}\right)  \right) \\  
&& = \frac{1}{s^2}+ \frac{f''\left(\frac{y-\beta^Tx}{s}\right)f\left(\frac{y-\beta^T x}{s}\right)-\left(f'\left(\frac{y-\beta^T x}{s}\right)\right)^2}{f\left(\frac{y-\beta^T x}{s}\right)^2}\cdot\left(\frac{y-\beta^Tx}{s^2}\right)^2  \\
&&\ + \  2 \ \frac{f'\left(\frac{y-\beta^T x}{s}\right)}{f\left(\frac{y-\beta^T x}{s}\right)}\cdot \left(\frac{y-\beta^T x}{s^3}\right). 
\end{eqnarray*}
Using similar arguments as above, we get
\begin{eqnarray*}
  -\mathbb E_{\theta_0} \left[ \frac{\partial^2 \log f_s(Y - \beta^T X)}{\partial s^2}|_{\theta = \theta_0} \right]   =  \frac{1}{s_0^2} \left(\int \frac{(f'(t))^2}{f(t)} t^2 dt  -1 \right) = c_2
\end{eqnarray*}
which is finite by (A11).  We conclude from the calculations above that $M_0 = I_0$.  Using similar calculations, we can easily show that
\begin{eqnarray*}
\left(
\begin{array}{c}
\mathbb E_{\theta_0}\left( \frac{\partial \log f_s(Y-\beta^T X)}{\partial \beta}|_{\theta = \theta_0} \right) \\
\mathbb E_{\theta_0}\left( \frac{\partial \log f_s(Y-\beta^T X)}{\partial s}|_{\theta = \theta_0} \right)
\end{array}
\right) = \mathbf{O}_{d+1}
\end{eqnarray*}
and 
\begin{eqnarray*}
\hspace{-1cm} && \left(
\begin{array}{cc}
\mathbb E_{\theta_0}\left[\left( \frac{\partial \log f_s(Y-\beta^T X)}{\partial \beta}|_{\theta = \theta_0} \frac{\partial \log f_s(Y-\beta^T X)}{\partial \beta^T}|_{\theta = \theta_0}\right)\right] &   \mathbb E_{\theta_0}\left[\left( \frac{\partial \log f_s(Y-\beta^T X)}{\partial \beta}|_{\theta = \theta_0} \frac{\partial \log f_s(Y-\beta^T X)}{\partial s}|_{\theta = \theta_0}\right)\right] \\
\mathbb E_{\theta_0}\left[\left( \frac{\partial \log f_s(Y-\beta^T X)}{\partial \beta^T}|_{\theta = \theta_0} \frac{\partial \log f_s(Y-\beta^T X)}{\partial s}|_{\theta = \theta_0}\right)\right] &  \mathbb E_{\theta_0}\left[\left( \frac{\partial \log f_s(Y-\beta^T X)}{\partial s}|_{\theta = \theta_0} \right)^2\right]
\end{array}
\right) \\
\hspace{-1cm}  && =    I_0.
\end{eqnarray*}
By applying the CLT we get 
\begin{eqnarray*}
Z_n:=\sqrt n \left( 
\begin{array}{c}
\frac{1}{n} \sum_{i=1}^n \frac{\partial \log f_s(Y_i - \beta^T X_i)}{\partial \beta}|_{\theta = \theta_0} \\
\frac{1}{n} \sum_{i=1}^n \frac{\partial \log f_s(Y_i - \beta^T X_i)}{\partial s}|_{\theta = \theta_0}  
\end{array}
\right) \to_d \mathcal N (\mathds{O}_{d+1}, I_0).
\end{eqnarray*}
From Assumption (A12) and the identity in (\ref{EqCTL}), it follows that
\begin{eqnarray*}
I_0^{-1} Z_n = -  \sqrt n (\widehat \theta_n - \theta_0) + o_{\mathbb P} (\sqrt n (\widehat \theta_n - \theta_0)).
\end{eqnarray*}
The result follows from Slutsky's Theorem.  \hfill $\Box$

\subsection{Proofs for Section \ref{sec: AsympEffic}}

\paragraph{Proof of the convergence in (\ref{asympOLS}).} 
We have that
\begin{eqnarray*}
\widehat{\beta}_{OLS}  = \left(\frac{1}{n}\sum_{i=1}^n X_i X_i^T\right)^{-1} \frac{1}{n}\sum_{i=1}^n X_i Y_i
\end{eqnarray*}
and hence
\begin{align*}
    \sqrt{n}(\widehat{\beta}_{OLS}-\beta_0) &= \sqrt{n}\left(\left(\sum_{i=1}^nX_iX_i^T\right)^{-1}\sum_{i=1}^nX_i(\beta_0^TX_i+ \epsilon_i) - \beta_0\right) \\
    & = \sqrt{n}\left(\left(\sum_{i=1}^nX_iX_i^T\right)^{-1}\sum_{i=1}^nX_i(X_i^T\beta_0+ \epsilon_i) - \beta_0\right) \\
    &= \sqrt{n}\left(\sum_{i=1}^nX_iX_i^T\right)^{-1}\sum_{i=1}^nX_i\epsilon_i\\
&=\left(\frac{1}{n}\sum_{i=1}^nX_iX_i^T\right)^{-1}\left(\frac{1}{\sqrt{n}}\sum_{i=1}^nX_i\epsilon_i\right). 
\end{align*}
By the Law of Large Numbers and the continuous mapping theorem, it holds that
$$
\left(\frac{1}{n}\sum_{i=1}^nX_iX_i^T\right)^{-1} \xrightarrow[]{\mathbb{P}} [\mathbb{E}[XX^T]]^{-1}.
$$
Note that $X_i \epsilon_i, i =1, \ldots, n$ are i.i.d. $d$-dimensional vectors such that $\mathbb E[X_1 \epsilon_1] = \mathbb E(X_1) \mathbb E(\epsilon_1) = 0$ using independence of $X$ and $\epsilon$. 
\begin{eqnarray*}
    V[X_1 \epsilon_1] &=&  \mathbb{E}[\epsilon^2XX^T] \\
    &=& \mathbb{E}[\epsilon^2]\mathbb{E}[XX^T] = \sigma^2\mathbb{E}[XX^T].  
\end{eqnarray*}
The claimed weak convergence follows from the Central Limit and Slutsky's theorems. 

\paragraph{Proof of Proposition \ref{asybehavols}.}  Using the definition of the asymptotic relative efficiency as given in \cite{serfling2009}, it follows from Theorem \ref{asympMLE} and the weak convergence in (\ref{asympOLS}) that
\begin{eqnarray*}
\eta & = &   \frac{c^{-1}_1}{\sigma^2}  = \frac{s^2_0}{\sigma^2 \int (f'(t))^2/f(t) dt}.
\end{eqnarray*}
Since $\sigma^2 = \text{var}(\epsilon)$, we must have that
\begin{eqnarray*}
\sigma^2 =  \int u^2 \frac{1}{s_0} f\left( \frac{u}{s_0} \right) du =  s_0^2  \int t^2 f(t) dt, \ \text{using $t = u/s_0$}
\end{eqnarray*}
which shows that $\eta$ is given by the claimed expression. Next, we will show that
\begin{eqnarray*}
\left(\int t^2 f(t) dt \right) \left(\int \frac{(f'(t))^2}{f(t)} dt\right) \ge 1
\end{eqnarray*}
with equality if and only if $f$ is the density of a centered Gaussian. Using the Cauchy Schwarz inequality we can write
\begin{eqnarray*}
\left(\int t f'(t) dt\right)^2    & =  &   \left(\int t \frac{f'(t)}{f(t)}  f(t) dt\right)^2   \\
& \le &  \left(\int t^2 f(t) dt \right)  \left( \int \left(\frac{f'(t)}{f(t)}\right)^2 f(t) dt\right) \\
& =  & \left(\int t^2 f(t) dt \right)  \left( \int \frac{(f'(t))^2}{f(t)} dt  \right), 
\end{eqnarray*}
and $\int t f'(t) dt =  [t f(t)]_{-\infty}^\infty - \int f(t) dt = -1$, using the fact that $\int \vert t \vert f(t) dt < \infty$.  This shows the claimed inequality.  Equality occurs if and only if there exists $\lambda \in \mathbb R \setminus \{0 \}$ such that for a.e. $t$
\begin{eqnarray*}
\frac{f'(t)}{f(t)}  = \lambda t,
\end{eqnarray*}
which holds if and only if $f(t)  =  \alpha \exp(\lambda t^2/2) $ for some $\alpha \ge 0$. Hence, $\lambda < 0$ and $f$ is the density of a centered Gaussian (with variance $2 \vert \lambda \vert$.  \hfill $\Box$

\paragraph{Derivation of $\eta_i(\gamma)$ for $i=1, 2, 3$.}
For the three families, we have that
$$\eta_i(\gamma)=\left(\int_\mathbb{R}\frac{(f'_{\gamma,i}(t))^2}{f_{\gamma,i}(t)}dt\right)^{-1}$$ since $\int_\mathbb{R}t^2f_{\gamma,i}(t)dt=1$, $\forall i \in \{1,2,3\}$.\\

\par \noindent \textbf{Family 1.}  For this family, the densities are given by
$$f_{\gamma,1}(t)=d_\gamma\exp(-c_\gamma|t|^\gamma),  t \in \mathbb{R}.$$
To compute $\eta_1(\gamma)$, we need to find the expression of $$\int_\mathbb{R}\frac{(f'_{\gamma,1}(t))^2}{f_{\gamma,1}(t)}dt.$$
By symmetry,
$$\int_\mathbb{R}\frac{(f'_{\gamma,1}(t))^2}{f_{\gamma,1}(t)}dt=2\int_0^{+\infty}\frac{(f'_{\gamma,1}(t))^2}{f_{\gamma,1}(t)}dt$$
For $t>0$, $$f_{\gamma,1}(t)=d_\gamma\exp(-c_\gamma t^\gamma).$$
Thus, $$f'_{\gamma,1}(t)=-d_\gamma c_\gamma \gamma t^{\gamma-1} \exp(-c_\gamma t^\gamma)$$
We have,
\begin{eqnarray*}
    \int_0^{+\infty}\frac{(f'_{\gamma,1}(t))^2}{f_{\gamma,1}(t)}dt &=& \int_0^{+\infty} d_\gamma c_\gamma^2\gamma^2t^{2\gamma-2}\exp(-c_\gamma t^\gamma)dt
\end{eqnarray*}
By the change-of-variable $x=t^\gamma \Leftrightarrow t=x^{\frac{1}{\gamma}}$, and $dt=\frac{1}{\gamma}x^{\frac{1}{\gamma}-1}dx$. Thus,
\begin{eqnarray*}
    \int_0^{+\infty}\frac{(f'_{\gamma,1}(t))^2}{f_{\gamma,1}(t)}dt &=& \int_0^{+\infty} d_\gamma c_\gamma^2\gamma^2x^{\frac{1}{\gamma}(2\gamma-2)}\exp(-c_\gamma x)\frac{1}{\gamma}x^{\frac{1}{\gamma}-1}dx \\
    &=& d_\gamma c_\gamma^2\gamma \int_0^{+\infty} x^{2-\frac{2}{\gamma}+\frac{1}{\gamma}-1}\exp(-c_\gamma x)dx \\
    &=& d_\gamma c_\gamma^2\gamma \int_0^{+\infty} x^{2-\frac{1}{\gamma}-1}\exp(-c_\gamma x)dx \\
    &=& d_\gamma c_\gamma^2\gamma \frac{\Gamma\left(2-\frac{1}{\gamma}\right)}{c_\gamma^{2-\frac{1}{\gamma}}}, \  \   \gamma \ \text{such that} \ 2-\frac{1}{\gamma}\geq 0 \\
    \Leftrightarrow \int_0^{+\infty}\frac{(f'_{\gamma,1}(t))^2}{f_{\gamma,1}(t)}dt &=& d_\gamma c_\gamma^\frac{1}{\gamma}\gamma \Gamma\left(2-\frac{1}{\gamma}\right), \  \   \gamma \in (2,\infty)
\end{eqnarray*}

We know that $d_\gamma=\frac{\gamma}{2}\left(\frac{\Gamma\left(\frac{3}{\gamma}\right)}{\Gamma\left(\frac{1}{\gamma}\right)^3}\right)^\frac{1}{2}$ and $c_\gamma=\left(\frac{\Gamma\left(\frac{3}{\gamma}\right)}{\Gamma\left(\frac{1}{\gamma}\right)}\right)^\frac{\gamma}{2}$

Then,
\begin{eqnarray*}
    \int_\mathbb{R}\frac{(f'_{\gamma,1}(t))^2}{f_{\gamma,1}(t)}dt&=&2d_\gamma c_\gamma^\frac{1}{\gamma}\gamma \Gamma\left(2-\frac{1}{\gamma}\right), \  \   \gamma \in \left(\frac{1}{2},\infty\right) \\
    &=& \gamma^2 \left(\frac{\Gamma\left(\frac{3}{\gamma}\right)}{\Gamma\left(\frac{1}{\gamma}\right)}\right)^\frac{1}{2}\left(\frac{\Gamma\left(\frac{3}{\gamma}\right)}{\Gamma\left(\frac{1}{\gamma}\right)^3}\right)^\frac{1}{2}\Gamma\left(2-\frac{1}{\gamma}\right) \\
    &=& \gamma^2 \frac{\Gamma\left(\frac{3}{\gamma}\right)}{\Gamma\left(\frac{1}{\gamma}\right)^2}\Gamma\left(2-\frac{1}{\gamma}\right), \  \  \gamma \in \left(\frac{1}{2},\infty\right) \\
    \Rightarrow \eta_1(\gamma)&=&\frac{\Gamma\left(\frac{1}{\gamma}\right)^2}{\gamma^2 \ \Gamma\left(\frac{3}{\gamma}\right) \Gamma\left(2-\frac{1}{\gamma}\right)}, \  \  \gamma \in (1/2, \infty)
\end{eqnarray*}

\par \noindent \textbf{Family 2.}
The densities are given by 
$$f_{\gamma,2}(t)   =  \displaystyle d_\gamma  \vert t \vert^{\gamma -1} \exp(- c_\gamma \vert t \vert), \gamma \ge 1.$$
For $t>0$,
\begin{eqnarray*}
    f_{\gamma,2}(t)   &=&  \displaystyle d_\gamma  t^{\gamma -1} \exp(- c_\gamma t) \\
    f'_{\gamma,2}(t) &=& \displaystyle d_\gamma ((\gamma-1)t^{\gamma-2}\exp(- c_\gamma t)-c_\gamma t^{\gamma-1}\exp(- c_\gamma t))\\
    &=& \displaystyle d_\gamma t^{\gamma-2}\exp(- c_\gamma t)(\gamma-1-c_\gamma t)
\end{eqnarray*}
Thus,
\begin{eqnarray*}
    \int_\mathbb{R}\frac{(f'_{\gamma,2}(t))^2}{f_{\gamma,2}(t)}dt&=&2\int_0^{+\infty}\frac{(f'_{\gamma,2}(t))^2}{f_{\gamma,2}(t)}dt \\
    &=& 2\int_0^{+\infty}\frac{d_\gamma^2 t^{2\gamma-4}\exp(-2 c_\gamma t)(\gamma-1-c_\gamma t)^2}{d_\gamma  t^{\gamma -1} \exp(- c_\gamma t)}dt \\
    &=& 2d_\gamma \int_0^{+\infty} t^{\gamma-2-1}\exp(- c_\gamma t)((\gamma-1)^2-2(\gamma-1)c_\gamma t+c_\gamma^2t^2) dt\\
    &=& 2d_\gamma \Big[ (\gamma-1)^2\int_0^{+\infty} t^{\gamma-2-1}\exp(- c_\gamma t)dt-2(\gamma-1)c_\gamma\int_0^{+\infty}t^{\gamma+1-1}\exp(- c_\gamma t)dt \\
    && \ + \ c_\gamma^2\int_0^{+\infty}t^{\gamma-1}\exp(- c_\gamma t)dt\Big]\\
    &=& 2d_\gamma \left[ (\gamma-1)^2\frac{\Gamma(\gamma-2)}{c_\gamma^{\gamma-2}}-2(\gamma-1)\frac{c_\gamma \Gamma(\gamma-1)}{c_\gamma^{\gamma-1}}+c_\gamma^2\frac{\Gamma(\gamma)}{c_\gamma^\gamma}\right] \\
    &=& \frac{2d_\gamma}{c_\gamma^{\gamma-2}}\left((\gamma-1)^2\Gamma(\gamma-2)-2(\gamma-1)\Gamma(\gamma-1)+\Gamma(\gamma)\right), \ \ \gamma \in (2,\infty)
\end{eqnarray*}
Using the well-known property $\Gamma(a+1)=a\Gamma(a), \forall a >0$, we can write,
\begin{eqnarray*}
    \int_\mathbb{R}\frac{(f'_{\gamma,2}(t))^2}{f_{\gamma,2}(t)}dt&=& \frac{2d_\gamma}{c_\gamma^{\gamma-2}}\left(\frac{(\gamma-1)^2\Gamma(\gamma)}{(\gamma-1)(\gamma-2)}-2\Gamma(\gamma)+\Gamma(\gamma)\right)\\
    &=& \frac{2d_\gamma}{c_\gamma^{\gamma-2}}\left(\frac{\gamma-1}{\gamma-2}-1\right)\Gamma(\gamma)\\
    &=& \frac{2d_\gamma}{c_\gamma^{\gamma-2}}\frac{1}{\gamma-2}\Gamma(\gamma), \ \ \gamma \in (2,\infty)
\end{eqnarray*}
We replace 
\begin{eqnarray*}
c_\gamma =  \left(\frac{\Gamma(\gamma +2)}{\Gamma(\gamma)}\right)^{\frac{1}{2}}, \ \ \text{and} \ \  d_\gamma =   \frac{1}{2 \Gamma(\gamma)} \left(\frac{\Gamma(\gamma+2)}{\Gamma(\gamma)}\right)^{\frac{\gamma}{2}}.  
\end{eqnarray*}
to obtain,
\begin{eqnarray*}
    \int_\mathbb{R}\frac{(f'_{\gamma,2}(t))^2}{f_{\gamma,2}(t)}dt&=&\left(\frac{\Gamma(\gamma+2)}{\Gamma(\gamma)}\right)^{\frac{\gamma}{2}}\left(\frac{\Gamma(\gamma)}{\Gamma(\gamma+2)}\right)^{\frac{\gamma-2}{2}}\frac{1}{\gamma-2}\\
    &=& \frac{\Gamma(\gamma+2)}{\Gamma(\gamma)(\gamma-2)}\\
    \Rightarrow \eta_2(\gamma) &=&  \frac{\Gamma(\gamma) (\gamma -2)}{\Gamma(\gamma +2)}, \  \  \gamma \in (2, \infty)
\end{eqnarray*}

\par \noindent \textbf{Family 3.}  In this case, we have that
$$f_{\gamma,3}(t)  = \displaystyle   d_\gamma \vert t \vert^{\gamma -1}  \exp(- c_\gamma \vert t \vert^\gamma), \gamma \ge 1.$$
For $t>0$,
\begin{eqnarray*}
f_{\gamma,3}(t)  &=& \displaystyle   d_\gamma t^{\gamma -1}  \exp(- c_\gamma t^\gamma)\\
f'_{\gamma,3}(t)  &=& \displaystyle   d_\gamma((\gamma-1) t^{\gamma -2}  \exp(- c_\gamma t^\gamma)-c_\gamma \gamma t^{\gamma -1} t^{\gamma -1} \exp(- c_\gamma t^\gamma) )\\
&=& \displaystyle   d_\gamma t^{\gamma -2}  \exp(- c_\gamma t^\gamma) ((\gamma-1)-c_\gamma \gamma t^\gamma)\\
\end{eqnarray*}
Thus, and again by symmetry,
\begin{eqnarray*}
    \int_\mathbb{R}\frac{(f'_{\gamma,3}(t))^2}{f_{\gamma,3}(t)}dt&=&2d_\gamma\int_0^{+\infty}\frac{t^{2\gamma -4}  \exp(- 2c_\gamma t^\gamma)(\gamma-1-c_\gamma \gamma t^\gamma)^2}{t^{\gamma -1} \exp(- c_\gamma t^\gamma)}dt\\
    &=& 2d_\gamma\int_0^{+\infty}t^{\gamma -3}  \exp(- c_\gamma t^\gamma)(\gamma-1-c_\gamma \gamma t^\gamma)^2dt
\end{eqnarray*}
By the change-of-variable $x=t^\gamma \Leftrightarrow t=x^{\frac{1}{\gamma}}$, and $dt=\frac{1}{\gamma}x^{\frac{1}{\gamma}-1}dx$. Thus,
\begin{eqnarray*}
    \int_\mathbb{R}\frac{(f'_{\gamma,3}(t))^2}{f_{\gamma,3}(t)}dt&=&\frac{2d_\gamma}{\gamma}\int_0^{+\infty}x^{\frac{\gamma-3}{\gamma}}\exp(- c_\gamma x)x^{\frac{1}{\gamma}-1}(\gamma-1-c_\gamma \gamma x)^2dx\\
    &=& \frac{2d_\gamma}{\gamma}\int_0^{+\infty}x^{-\frac{2}{\gamma}}\exp(- c_\gamma x)((\gamma-1)^2-2c_\gamma \gamma(\gamma-1)x+c_\gamma^2\gamma^2x^2)dx\\
    &=& \frac{2d_\gamma}{\gamma}\left((\gamma-1)^2\int_0^{+\infty}x^{-\frac{2}{\gamma}}\exp(- c_\gamma x)dx\right.\\
    && \ - \ \left. 2c_\gamma \gamma(\gamma-1)\int_0^{+\infty}x^{1-\frac{2}{\gamma}}\exp(- c_\gamma x)dx + c_\gamma^2\gamma^2\int_0^{+\infty}x^{2-\frac{2}{\gamma}}\exp(- c_\gamma x)dx\right)\\
    &=&  \frac{2d_\gamma}{\gamma}\left((\gamma-1)^2\frac{\Gamma\left(-\frac{2}{\gamma} +1 \right)}{c_\gamma^{-\frac{2}{\gamma} +1}}-2c_\gamma \gamma(\gamma-1)\frac{\Gamma\left(2-\frac{2}{\gamma}\right)}{c_\gamma^{2-\frac{2}{\gamma} }}+ c_\gamma^2\gamma^2\frac{\Gamma\left(3-\frac{2}{\gamma}\right)}{c_\gamma^{3-\frac{2}{\gamma} }}\right)
\end{eqnarray*}
for , $\gamma \in (2,\infty)$. Using that,
\begin{eqnarray*}
    \Gamma\left(2-\frac{2}{\gamma}\right)= \left(1-\frac{2}{\gamma}\right)\Gamma\left(1-\frac{2}{\gamma}\right) \ \ \text{and} \ \ 
    \Gamma\left(3-\frac{2}{\gamma}\right)= \left(2-\frac{2}{\gamma}\right)\left(1-\frac{2}{\gamma}\right)\Gamma\left(1-\frac{2}{\gamma}\right)
\end{eqnarray*}
it follows that,
\begin{eqnarray*}
    \int_\mathbb{R}\frac{(f'_{\gamma,3}(t))^2}{f_{\gamma,3}(t)}dt&=&\frac{2d_\gamma}{\gamma}\frac{1}{c_\gamma^{-\frac{2}{\gamma}+1}}\left((\gamma-1)^2-2\gamma(\gamma-1)\left(1-\frac{2}{\gamma}\right)+\gamma^2\left(2-\frac{2}{\gamma}\right)\left(1-\frac{2}{\gamma}\right)\right)\Gamma\left(1-\frac{2}{\gamma}\right)
\end{eqnarray*}
We replace 
\begin{eqnarray*}
c_\gamma  =  \Gamma\left(\frac{2}{\gamma} +1 \right)^{\frac{\gamma}{2}}, \ \ \text{and} \ \  d_\gamma =   \frac{\gamma}{2}  \Gamma\left(\frac{2}{\gamma} +1 \right)^{\frac{\gamma}{2}}.  
\end{eqnarray*}
Finally,
\begin{eqnarray*}
    \int_\mathbb{R}\frac{(f'_{\gamma,3}(t))^2}{f_{\gamma,3}(t)}dt&=&\frac{\Gamma\left(\frac{2}{\gamma} +1 \right)^{\frac{\gamma}{2}}}{\Gamma\left(\frac{2}{\gamma} +1 \right)^{\left(1-\frac{2}{\gamma}\right)\frac{\gamma}{2}}}\left((\gamma-1)^2-2(\gamma-1)(\gamma-2)+2(\gamma-1)(\gamma-2)\right)\Gamma\left(1-\frac{2}{\gamma}\right)\\
    &=& \Gamma\left(\frac{2}{\gamma} +1 \right)^{\frac{\gamma}{2}-\frac{\gamma}{2}+1}(\gamma-1)^2\Gamma\left(1-\frac{2}{\gamma}\right)\\
    &=& \Gamma\left(\frac{2}{\gamma} +1 \right)\Gamma\left(1-\frac{2}{\gamma}\right)(\gamma-1)^2\\
    \Rightarrow \eta_3(\gamma) &=& \frac{1}{\Gamma\left(\frac{2}{\gamma} +1 \right) \Gamma\left(1-\frac{2}{\gamma}\right) (\gamma -1)^2} , \  \  \gamma \in (2, \infty).
\end{eqnarray*}



\end{document}